\documentclass[11pt]{article} 

\usepackage[utf8]{inputenc}

\usepackage{geometry}
\usepackage[dvips]{graphicx}

\usepackage{booktabs} 
\usepackage{array} 
\usepackage{paralist} 
\usepackage{verbatim} 
\usepackage{subfig} 
\usepackage{cite}

\usepackage{fancyhdr} 
\usepackage[nottoc,notlof,notlot]{tocbibind}
\usepackage[titles,subfigure]{tocloft}

\usepackage{amsmath,amsfonts,amsthm,mathrsfs,amssymb}
\usepackage{latexsym}
\usepackage{enumerate}
\usepackage[usenames]{color}

\newtheorem{thm}{Theorem}[section]

\newtheorem{lem}{Lemma}[section]
\newtheorem{prop}{Proposition}[section]
\theoremstyle{definition}
\newtheorem{defin}{Definition}[section]
\theoremstyle{remark}

\numberwithin{equation}{section}

\def\bb{\begin{equation}} \def\ee{\end{equation}}
\def\beqn{\begin{eqnarray}}  \def\eqn{\end{eqnarray}}
\def\beqnx{\begin{eqnarray*}} \def\eqnx{\end{eqnarray*}}

{

}

\newcommand{\R}{{\mathbb R}}
\def\ba{\begin {eqnarray*} }
\def\ea{\end {eqnarray*} }
\def\beq{\begin {eqnarray}}
\def\eeq{\end {eqnarray}}
\def\p{\partial}

\title{Regularized Transformation-Optics Cloaking in Acoustic and Electromagnetic Scattering}

\author{Hongyu Liu\thanks{Department of Mathematics, Hong Kong Baptist University, Kowloon, Hong Kong.   Email:  {\tt hongyu.liuip@gmail.com}},\ \ Gunther Uhlmann\thanks{Department of Mathematics, University of Washington, Seattle, WA 98195, USA and Fondation des Sciences Math\'ematiques de Paris.  Email:  {\tt gunther@math.washington.edu}}    }

\begin{document}

\date{}

\maketitle

\begin{abstract}
We consider transformation-optics based cloaking in acoustic and electromagnetic scattering. The blueprints for an ideal cloak use singular acoustic and electromagnetic materials, posing server difficulties to both theoretical analysis and practical fabrication. In order to avoid the singular structures, various regularized approximate cloaking schemes have been developed. We survey these developments in this paper. We also propose some challenging issues for further investigation.
\end{abstract}

\section{Introduction}

We shall be concerned with invisibility cloaking for acoustic and electromagnetic (EM) waves. A region is said to be {\it cloaked} if its content together with the cloak is indistinguishable from the background space with respect to exterior wave measurements. A proposal for cloaking for electrostatics using the invariance properties
of the conductivity equation was pioneered in \cite{GLU,GLU2}.
Blueprints for making objects invisible to electromagnetic (EM) waves were proposed in two articles in {\it Science} in 2006 \cite{Leo,PenSchSmi}. The article by Pendry et al
uses the same transformation used in \cite{GLU,GLU2} while the work of Leonhardt uses a conformal mapping in two dimensions. The method based on the invariance properties of the equations modelling the wave phenomenon has been named
{\sl transformation optics} and has received a lot of attention in the scientific community and the popular press because of the generality of the method and its simplicity. There have been several other proposals for cloaking. We mention the works of Milton and Nicorovici \cite{MN} and of Alu and Engetha \cite{AE}. 

The method of transformation optics relies on the transformation properties of optical parameters and the transformation invariance of the governing wave equations. To obtain an ideal invisibility cloak, one first selects a region $\Omega$ in the background space for constructing the cloaking device. Throughout the paper, we assume that the background space is uniformly homogeneous in order to facilitate the exposition, but all of the results discussed in this paper can be straightforwardly extended to the case with an inhomogeneous background space. Let $P\in\Omega$ be a single point and let $F$ be a diffeomorphism which blows up $P$ to a region $D$ within $\Omega$. The ambient homogeneous medium around $P$ is then `compressed' via the {\it push-forward} of the transformation to form the cloaking medium in $\Omega\backslash\overline{D}$, whereas the `hole' $D$ is the cloaked region within which one can place the target object. The cloaking region $\Omega\backslash\overline{D}$ and the cloaked region $D$ form the cloaking device in the {\it physical space}, whereas the homogeneous background space containing the singular point $P$ is referred to as the {\it virtual space} (see Figure 3). Due to the transformation invariance of the corresponding wave equations, the acoustic/EM scattering in the physical space with respect to the cloaking device is the same as that in the virtual space. Heuristically speaking, the scattering information of the cloaking device is then `hidden' in the singular point $P$. In a similar fashion cloaking devices based on blowing up a {\it crack} (namely, a curve in $\mathbb{R}^3$) or a {\it screen} (namely, a flat surface in $\mathbb{R}^3$) were, respectively, considered in \cite{GKLU2} and \cite{LiP}, resulting in the so-called EM wormholes and carpet-cloak respectively. 

The blow-up-a-point (respectively, -crack or -screen) construction yields singular cloaking materials, namely, the material parameters violate the regular conditions. The singular media present a great challenge for both theoretical analysis and practical fabrications. In order to tackle the acoustic and electromagnetic wave equations with singular coefficients underlying the ideal invisibility cloaks, finite energy solutions on Sobolev spaces with singular weights
were introduced and studied in \cite{GKLU3,GKLU2,HetLiu,LZ1}. On the other hand,several regularized constructions have been developed in the literature in order to avoid the singular structures. In \cite{GKLUoe,GKLU_2,RYNQ}, a truncation of singularities has been introduced. In \cite{KOVW,KSVW,Liu}, the blow-up-a-point transformation in \cite{GLU2,Leo,PenSchSmi} has been regularized to become the `blow-up-a-small-region' transformation. Nevertheless, it is pointed out in \cite{KocLiuSun} that the truncation-of-singularities construction and the blow-up-a-small-region construction are equivalent to each other.  Instead of ideal/perfect invisibility, one would consider approximate/near invisibility for a regularized construction; that is, one intends to make the corresponding wave scattering effect due to a regularized cloaking device as small as possible depending on an asymptotically small regularization parameter $\rho\in\mathbb{R}_+$.

Due to its practical importance, the approximate cloaking has recently been extensively studied. In \cite{Ammari1,KSVW}, approximate cloaking schemes were developed for electrostatics. In \cite{Ammari2,Ammari3,KOVW,LiLiuRonUhl,LiLiuSun,Liu,Liu2,LiuSun,N1,N2}, various near-cloaking schemes were presented for scalar waves governed by the Helmholtz equation. In \cite{Ammari4,BL,BLZ,LiuZhou}, near-cloaking schemes were developed and investigated for the vector waves governed by the Maxwell system. Generally speaking, a regularized near-cloak consist of three layers: the innermost core is the cloaked region, the outermost layer is the cloaking region, and a compatible lossy layer right between the cloaked and cloaking regions. In the cloaking layer, the cloaking parameters are obtained by the
push-forward construction mentioned earlier. Inside the cloaked region, from a practical viewpoint, one can place an arbitrary content, which could be both passive and active. The special lossy layer employed right between the cloaked and cloaking regions has shown to be necessary \cite{KOVW,LiuZhou}, since otherwise there exist cloak-busting inclusions which defy any attempt for cloaking at particular resonant frequencies. In the extreme case when the lossy parameters go to infinity, the lossy layer become an impenetrable obstacle layer, and this is the one considered in \cite{Ammari2,Ammari3,Ammari4,Liu}. In the rest of this paper, we shall survey these developments and at certain places, we shall also point out challenges for  further investigation. In addition to the present survey, we also refer to the survey papers \cite{CC,GKLU4,GKLU5,U2,YYQ} and the references therein for discussions of the theoretical and experimental progress on invisibility cloaking. In this paper we make emphasis on remote observations via the scattering amplitude or scattering operator. The same considerations are valid for the ``near-field" which corresponds to the Cauchy data or the Dirichlet-to-Neumann map \cite{U}. In Section 2, we review perfect cloaking for the case of electrostatics using Cauchy data.
In Section 3 we define precisely what we mean by perfect cloaking for scattering. In Section 4 we review the push-forward construction. In Section 5 we consider regularized or approximate cloaks for acoustics, including the case of partial cloaks. In Section 6 we discuss regularized cloaks in electromagnetics.

\section{Invisibility
for electrostatics}\label{sec-to for es}
We discuss here only perfect cloaking for electrostatics. For
similar results for electromagnetic waves, acoustic waves, quantum
waves, etc.,  see the review papers \cite{GKLU2}, \cite{GKLU4} and the references
given there.
The fact that  the boundary measurements do not change, when a
conductivity is  pushed forward by a smooth diffeomorphism  leaving
the boundary fixed can already be considered as a
weak form of invisibility. Different media appear to be the same,
and the apparent location of objects  can change. However, this does
not yet constitute real invisibility, as nothing has been hidden
from view.
In  invisibility cloaking the aim is to hide an object inside a
domain by surrounding it with a material so that even the
presence of this object can not be detected by measurements on the
domain's boundary. This means that
 all boundary
measurements for the domain with this cloaked object included would
be the same as if the domain were filled with a homogeneous,
isotropic material. Theoretical models for this have been found by
applying diffeomorphisms having singularities. These were first
introduced in the framework of electrostatics, yielding
counterexamples to the anisotropic Calder\'on problem (see \cite{U1} for a review) in the form of
singular, anisotropic conductivities in $\R^n, n\ge 3$,
indistinguishable from a constant isotropic conductivity in that
they have the same Dirichlet-to-Neumann map  \cite{GLU, GLU2}. The
same construction was rediscovered for electromagnetism in
\cite{PenSchSmi}, with the intention of actually building such a device
with appropriately designed metamaterials; a modified version of
this was then experimentally demonstrated in \cite{Sc}. (See also
\cite{Leo} for a somewhat different approach to cloaking in the high
frequency limit.)
The first constructions in this direction were based on  blowing up
the metric around a point \cite{LTU}. In this construction, let
$(M,g)$ be a compact 2-dimensional manifold with non-empty boundary,
let $x_0\in M$ and consider the manifold \ba \widetilde M=M\setminus
\{x_0\} \ea with the metric \ba \widetilde g_{ij}(x)=\frac
1{d_M(x,x_0)^2}g_{ij}(x), \ea where $d_M(x,x_0)$ is the distance
between $x$ and $x_0$ on $(M,g)$. Then $(\widetilde M,\widetilde g)$ is a
complete, non-compact 2-dimensional Riemannian manifold with the
boundary $\p \widetilde M=\p M$. Essentially, the point $x_0$ has been
``pulled to infinity''. On the manifolds $M$ and $\widetilde M$ we
consider the boundary value problems \ba \left\{\begin{array}{l}
\Delta_g u=0\quad \hbox{in $M$,}\\
u=f\quad \hbox{on $\p M$,}\end{array}\right. \quad\hbox{and}\quad
\left\{\begin{array}{l}
\Delta_{\widetilde g} \widetilde u=0\quad \hbox{in $\widetilde M$,}\\
\widetilde u=f\quad \hbox{on $\p \widetilde M$,}\\
\widetilde u\in L^\infty(\widetilde M).\end{array}\right. \ea  Here $\Delta_g$ denotes the Laplace-Beltrami operator associated to the metric $g$. Note that in dimension $n\ge 3$ metrics and conductivities are equivalent. These boundary
value problems are uniquely solvable and define the DN maps \ba
\Lambda_{M,g}f=\p_\nu u|_{\p M},\quad \Lambda_{\widetilde M,\widetilde
g}f=\p_\nu \widetilde u|_{\p \widetilde M} \ea where $\p_\nu$ denotes the
corresponding conormal derivatives. Since, in the two dimensional
case, functions which are harmonic with respect to the metric $g$
stay harmonic with respect to any
 metric which is conformal to $g$,
one can see that $\Lambda_{M,g}=\Lambda_{\widetilde M,\widetilde g}$. This
can be seen using e.g.\ Brownian motion or capacity arguments. Thus,
the boundary measurements for $(M,g)$ and $(\widetilde M,\widetilde g)$
coincide. This gives a counterexample for the inverse
 electrostatic problem on Riemannian surfaces -- even the topology
of possibly non-compact Riemannian surfaces can not be determined
using boundary measurements (see Fig.\ 1).
\begin{figure}[htbp]\label{LTU figure}
\begin{center}
\includegraphics[width=8cm]{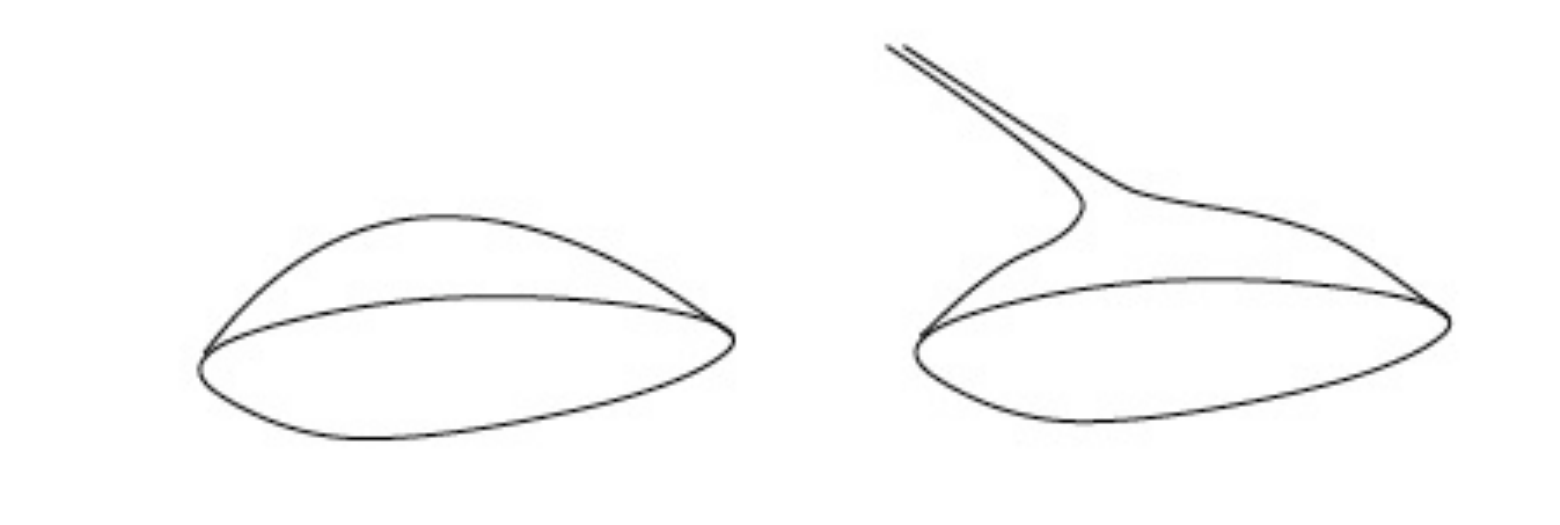} \label{two}
\caption{Blowing up a metric at a point, after \cite{LTU}. The
electrostatic boundary measurements on the boundary of the
surfaces, one compact and the other noncompact but complete,
coincide.}
\end{center}
\end{figure}
The above example can be thought as a ``hole'' in a Riemann surface
that does not change the boundary measurements. Roughly speaking,
mapping the manifold $\widetilde M$ smoothly to the set $M\setminus
\overline B_M(x_0,\rho)$, where $B_M(x_0,\rho)$ is a metric ball of
$M$, and by putting an object in the obtained hole $\overline
B_M(x_0,\rho)$, one could hide it from detection at the boundary.
This {observation was used} in \cite{GLU,GLU2}, where
``undetectability" results were  introduced in three dimensions,
using degenerations of Riemannian metrics, whose singular limits can
be considered as coming directly from singular changes of variables.
\begin{figure}[htbp]\label{collapce}
\begin{center}
\includegraphics[width=4cm]{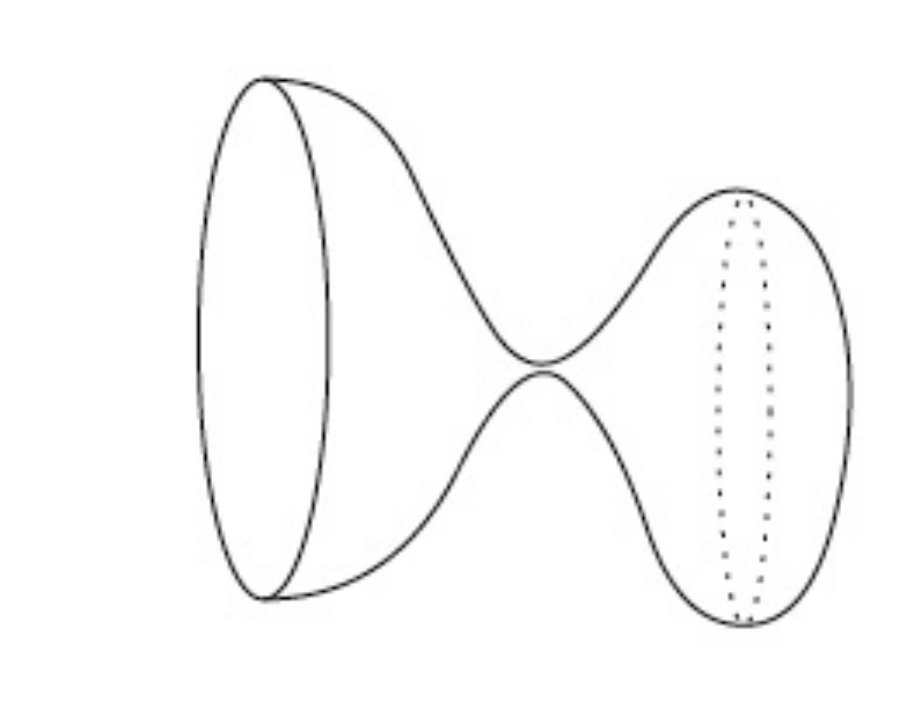}
\caption{A typical member of a family of manifolds  developing a
singularity as the width of the neck connecting the two parts goes
to zero.}
\end{center}
\end{figure}
The degeneration of the metric (see Fig. 2) can be obtained by
considering surfaces (or  manifolds in the higher dimensional cases)
with a thin ``neck'' that is pinched. At the limit the manifold
contains a pocket about which the boundary measurements do not give
any information. If the collapsing of the manifold is done in an
appropriate way, we have, in the limit, a singular Riemannian
manifold which is indistinguishable in boundary measurements  from a
flat surface. Then the conductivity which corresponds to this metric
is also singular at the pinched points, cf. the first formula in
(\ref{eq: cond}). The electrostatic measurements on the boundary for
this singular  conductivity will be the same
 as for  the original regular  conductivity
corresponding to the metric $g$.
To give a  precise, and concrete, realization of this idea, let
$B(0,R)\subset \R^3$ denote the open ball with center 0 and radius
$R$. We use in the sequel the set $N=B(0,2)$, the region at the
boundary of which the electrostatic measurements will be made,
decomposed into two parts, $N_1=B(0,2)\setminus \overline B(0,1)$
and  $N_2=B(0,1)$. We call the interface $\Sigma=\p N_2$  between
$N_1$ and $N_2$ the {\it cloaking surface}.
We  also use  a ``copy'' of the ball $B(0,2)$, with the notation
$M_1=B(0,2)$,  another ball $M_2=B(0,1)$, and the disjoint union $M$
of $M_1$ and $M_2$. (We will see the reason for distinguishing
between $N$ and $M$.) Let $g_{jk}=\delta_{jk}$ be the Euclidian
metrics  in $M_1$ and $M_2$ and let $\gamma=1$ be the corresponding
isotropic homogeneous conductivity. We define a singular
transformation \beq \label{eqn-sing transf} F_1:M_1\setminus\{0\}\to
N_1 ,\quad
   \ F_1(x)=(\frac {|x|}2+1)\frac x{|x|},\quad
0<|x|\le 2. \eeq
\begin{figure}[htbp]\label{figure for map F_1.}
\begin{center}
\includegraphics[width=.7\linewidth]{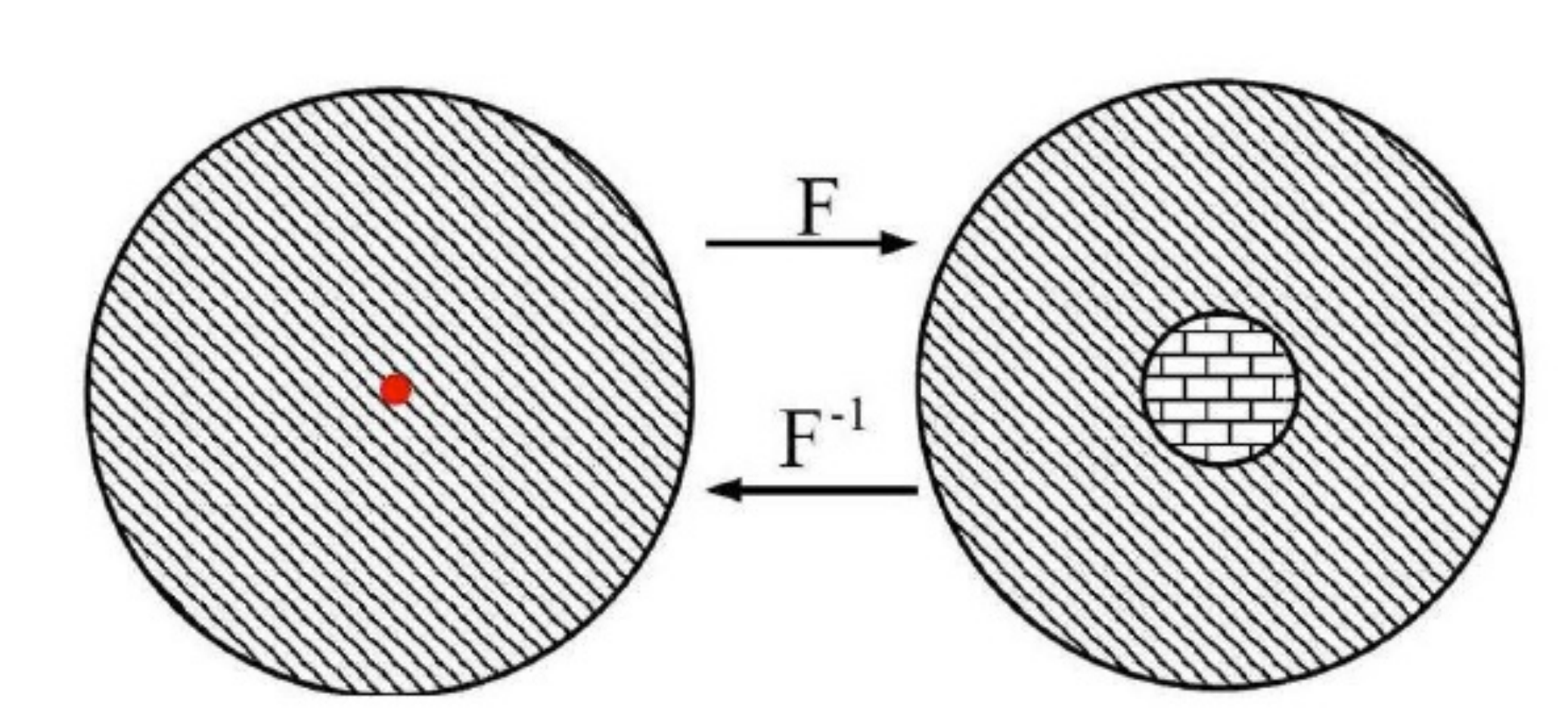}
\caption{Map $F_1:B(0,2)\setminus\{0\}\to B(0,2)\setminus\overline
B(0,1)$}\end{center}\end{figure}
We also consider a regular transformation (diffeomorphism) $F_2: M_2
\mapsto N_2$, which for simplicity we take to be the identity map
$F_2=Id$. Considering the maps $F_1$ and $F_2$ together, $F= (F_1,
F_2)$, we define a map $F:M\setminus \{0\}= (M_1\setminus \{0\})\cup
M_2 \to N\setminus \Sigma$.
The push-forward $\widetilde g=F_*g$ of the metric $g$ in $M$ by $F$ is
the metric in $N$ given by \beq\label{transf 2}
\left(F_*g\right)_{jk}(y)=\left. \sum_{p,q=1}^n \frac {\p F^p}{\p
x^j}(x) \,\frac {\p F^q}{\p x^k}(x)  g_{pq}(x)\right|_{x=F^{-1}(y)}.
\eeq This metric gives rise to a  conductivity $\widetilde \gamma$ in
$N$ which is singular in $N_1$, \beq\label{eq: cond} \widetilde \gamma
=\left\{\begin{array}{ll} |\widetilde g|^{1/2}\widetilde g^{jk}
  & \hbox{for }x\in N_1,\\
\delta^{jk} & \hbox{for }x\in N_2.\end{array}\right. \eeq
Thus, $F$ forms an invisibility construction that we call ``blowing
up a point". Denoting by $(r,\phi,\theta)\mapsto
      (r\sin\theta \cos \phi,r\sin\theta \sin \phi,r\cos\theta)$ the
spherical coordinates, we have
\begin{equation}\label{eqn-sing tensor 2}
\widetilde \gamma= \left(\begin{array}{ccc}
2(r-1)^2\sin \theta & 0 & 0\\
0 & 2 \sin \theta & 0 \\
0 & 0 &  2 (\sin \theta)^{-1}\\
\end{array}
\right), \quad 1<|x|\leq 2.
\end{equation}
Note that the anisotropic conductivity $\widetilde \gamma$ is singular
(degenerate) on $\Sigma$ in the sense that it is not bounded from
below by any positive multiple of $I$ (see \cite{KSVW} for a
similar calculation).
\begin{figure}[htbp]\label{Analytic sol.}
\begin{center}
\includegraphics[width=.5\linewidth]{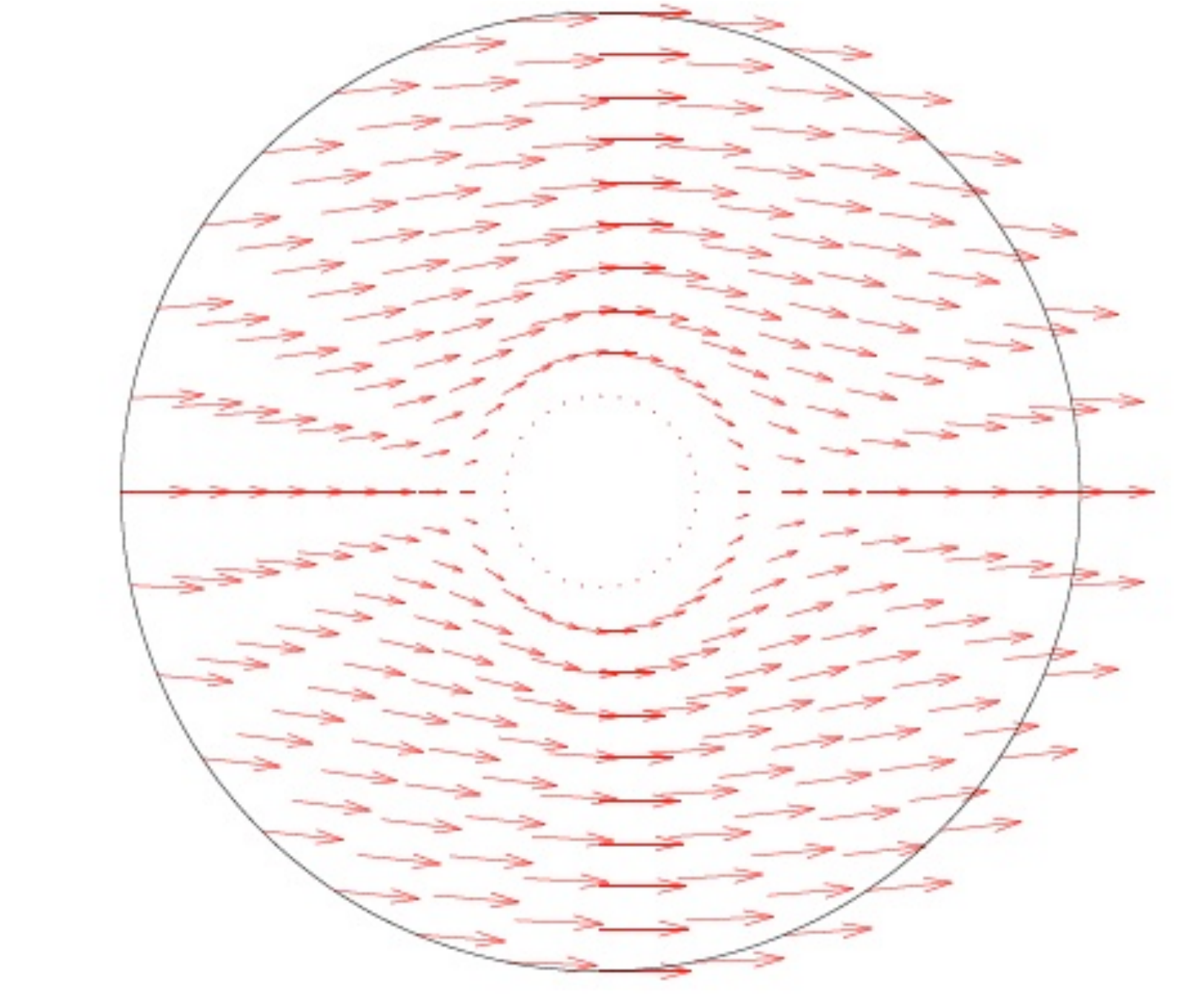}
\caption{Analytic solutions for the currents}
\end{center}
\end{figure}
The Euclidean conductivity $\delta^{jk}$  in $N_2$ (cf. \eqref{eq: cond})
could be replaced by any smooth conductivity bounded from below and
above by positive constants. This would correspond to cloaking of a
general object with  non-homogeneous, anisotropic conductivity.
Here, we use the Euclidian metric just for simplicity.
Consider now the {\emph{Cauchy data}} of all solutions in the
Sobolev space $H^1(N)$ of the conductivity equation corresponding to
$\widetilde \sigma$, that is, \ba C_1(\widetilde \gamma)=\{(u|_{\p
N},\nu\cdotp \widetilde \gamma \nabla u|_{\p N})\ :\ u\in H^1(N),\
\nabla\cdotp \widetilde \gamma\nabla u=0\}, \ea where $\nu$ is the
Euclidian unit normal vector of $\p N$.
\begin{thm}\label{thm-cond}
(\cite{GLU2}) The Cauchy data of all $H^1$-solutions for the
conductivities $\widetilde \gamma$ and $\gamma_0:=1$ on $N$ coincide, that
is, $C_1(\widetilde \gamma)=C_1(\gamma_0)$.
\end{thm}
This means that all boundary measurements for
  the homogeneous conductivity $\gamma_0=1$ and the degenerated
conductivity $\widetilde \gamma$  are the same. The result above was
proven in \cite{GLU,GLU2} for the case of dimension $n\ge 3.$ The
same basic construction works in the two dimensional case
\cite{KSVW}.
Fig.\ 4 portrays an analytically obtained solution on a disc with
conductivity $\widetilde \gamma$. As seen in the figure, no currents
appear near the center of the disc, so that if the conductivity is
changed near the center, the measurements on the boundary $\p N$ do
not change.
The above invisibility result  is valid for a more general class of
singular cloaking transformations. A general class, sufficing at
least for electrostatics, is given by the following result from
\cite{GLU2}:
\begin{thm}\label{main B}
Let $\Omega\subset \R^n$, $n\geq 3$, be a bounded smooth domain and $g=(g_{ij})$ be a smooth
metric on $\Omega$ bounded from above and below by positive
constants. Let $D\subset\subset \Omega$ be a smooth subdomain and such that there is a
$C^\infty$-diffeomorphism $F:\Omega\setminus\{y\}\to \Omega\setminus
\overline D$ satisfying $F|_{\p \Omega}=Id$ and such that
\beq\label{Q 4} dF(x)\geq c_0I,\quad \hbox{det}\,(DF(x))\geq
c_1\,\hbox{dist}_{{}_{\R^n}}\,(x,y)^{-1} \eeq where $DF$ is the
Jacobian matrix in Euclidian coordinates on $\R^n$ and $c_0,c_1>0$.
Let $\hat g$ be a metric in $\Omega$ which coincides with $\widetilde
g=F_*g$ in $\Omega\setminus \overline D$ and is an arbitrary regular
positive definite metric in $D^{int}$.
 Finally, let  $\gamma$ and $\hat \gamma$
be the conductivities corresponding to $g$ and $\hat g$ (cf.\
(\ref{eq: cond})). Then, \ba C_1(\hat \gamma)=C_1(\gamma). \ea
\end{thm}
The key to the proof of Theorem \ref{main B} is a removable
singularities theorem that implies that solutions of the
conductivity equation in  $\Omega\setminus \overline D$ pull back by
this singular transformation to solutions of the conductivity
equation in the whole $\Omega$.
Returning to the case $\Omega=N$ and the conductivity given by
(\ref{eq: cond}), similar types of results are valid also for a more
general class of solutions. Consider an unbounded quadratic form,
$A$ in $L^2(N, |\widetilde g|^{1/2} dx)$, \ba A_{\widetilde
\gamma}[u,v]=\int_N\widetilde\gamma\nabla u\cdotp \nabla v\,dx \ea
defined for $u,v\in {\mathcal D}(A_{\widetilde \gamma})=C_0^\infty(N)$.
Let $\overline A_{\widetilde \gamma}$ be the closure of this quadratic
form and say that
\ba \nabla\cdotp \widetilde \gamma\nabla u=0\quad\hbox{in }N \ea is
satisfied in the finite energy sense if there is $u_0\in H^1(N)$
supported in $N_1$ such that $u-u_0\in {\mathcal D}(\overline
A_{\widetilde \gamma})$ and \ba \overline A_{\widetilde \gamma}[u-u_0,v]=
-\int_N\widetilde\gamma\nabla u_0\cdotp \nabla v\,dx,\quad \hbox{for all
}v\in {\mathcal D}(\overline A_{\widetilde \gamma}). \ea Then the Cauchy
data set of the finite energy solutions, denoted by \ba
C_{f.e.}(\widetilde \gamma)=\Big\{(u|_{\p N},\nu\cdotp \widetilde \gamma
\nabla u|_{\p N})\: | \hbox{$u$ is a finite energy solution of
$\nabla\cdotp \widetilde \gamma\nabla u=0$}\Big\}, \ea coincides with
the Cauchy data $C_{f.e.}(\gamma_0)$ corresponding to the homogeneous
conductivity $\gamma_0=1$, that is, \beq\label{fin. en. solutions}
C_{f.e.}(\widetilde \gamma)=C_{f.e.}(\gamma_0). \eeq
Kohn, Shen, Vogelius and Weinstein have considered in \cite{KSVW}  the case when instead of blowing up a point
one stretches a small ball into the cloaked region (see Fig.~5). In this case the
conductivity is non-singular and one gets ``almost" invisibility
with a precise estimate in terms of the radius of the small ball. We study almost cloaking in more detail in the next sections for acoustic and electromagnetic scattering.

\section{Wave scattering and invisibility cloaking}

We start with the direct scattering for acoustic waves. Let $\Omega$ and $D$ be two bounded Lipschitz domains in $\mathbb{R}^n$, $n\geq 2$, such that $D\Subset \Omega$. 
Let $\eta=\eta(x)=(\eta^{ij}(x))\in\mathbb{R}^{n\times n}$, $x\in\mathbb{R}^n$, be a symmetric-matrix valued measurable function such that, for some $\lambda$, $0<\lambda\leq 1$, we have
\begin{equation}\label{eq:reg1}
\lambda |\xi|^2\leq\sum_{i,j=1}^n\eta^{ij}(x)\xi_i\xi_j\leq \lambda^{-1} |\xi|^2\quad \mbox{for any $\xi\in\mathbb{R}^n$ and for a.e. $x\in\mathbb{R}^n$}.
\end{equation}
Let $q=q_1+i q_2=q(x)$, $x\in\mathbb{R}^n$, be a complex-valued bounded measurable function with real and imaginary parts $q_1$ and $q_2$ respectively, such that, for some $\lambda$, $0<\lambda\leq 1$, we have
\begin{equation}\label{eq:reg2}
q_1(x)\geq \lambda,\quad q_2(x)\geq 0\quad \mbox{for a.e. $x\in\mathbb{R}^n$}.
\end{equation}
Furthermore, we assume that $q(x)=q_0:=1$ and $\eta^{ij}(x)=\eta_0^{ij}:=\delta^{ij}$ for $x\in\mathbb{R}^n\backslash\overline{\Omega}$. In the following, \eqref{eq:reg1} and \eqref{eq:reg2} will be referred to as the {\it regular conditions} on $\eta$ and $q$, and $\lambda$ is called the {\it regular constant}.

Next, we introduce the time-harmonic wave scattering governed by the Helmholtz equation
whose weak solution is $u=u(x,d,k)$, $x\in\mathbb{R}^n$, where $d\in\mathbb{S}^{n-1}$, $k\in\mathbb{R}_+$,
\begin{equation}\label{eq:Helmholtz II}
\begin{cases}
\ \displaystyle{\mathrm{div}(\eta\nabla u)+{k^2}q u}=-h & \hspace*{-1.7cm} \text{in }\mathbb{R}^n,\\
\ \text{$u(x,d,k)-e^{ik\cdot d}$ satisfies the radiation condition},
\end{cases}
\end{equation}
where $h$ is supported in $\Omega$ and $h\in L^2(\Omega)$.
The statement in \eqref{eq:Helmholtz II} means that if one lets $u^s(x,d,k)=u(x,d,k)-e^{ikx\cdot d}$, then
\begin{equation}\label{eq:sommerfeldn}
\lim_{r\rightarrow\infty}r^{\frac{n-1}{2}}\left(\frac{\partial u^s(x)}{\partial r}-i k u^s(x)\right )=0, \quad r=|x|.
\end{equation}

In the physical situation, \eqref{eq:Helmholtz II} can be used to describe the time-harmonic acoustic scattering due to an inhomogeneous acoustical medium $(\Omega; \eta, q)$ located in an otherwise uniformly homogeneous space $(\mathbb{R}^n\backslash\overline{\Omega}; \eta_0, q_0)$. $\eta^{-1}$ and $q$, respectively, denote the density tensor and acoustic modulus of the acoustical medium, and the RHS term $h$ denotes a source/sink in $\Omega$. $u(x)$ is the wave pressure with $U(x,t):=u(x) e^{-ik t}$ representing the wave field satisfying the scalar wave equation
\[
q(x)U_{tt}(x,t)-\sum_{i,j=1}^N \frac{\partial}{\partial x_i}\left(\eta^{ij}(x) \frac{\partial}{\partial x_j}U(x,t) \right )=h(x)e^{-ikt},\quad\mbox{\ \ $(x,t)\in \mathbb{R}^n\times \mathbb{R}$}.
\]
The function $u^i(x):=e^{ik x\cdot d}$ is an incident plane wave with $k$ denoting the wave number and $d\in\mathbb{S}^{n-1}$ denoting the impinging direction. $u(x)$ is called the total wave field and $u^s(x)$ is called the scattered wave field, which is the perturbation of the incident plane wave caused by the presence of the inhomogeneity $(\Omega; \eta, q)$ in the whole space. Indeed, it is easily seen that if there is no presence of the inhomogeneity, $u^s$ will be vanishing. 

We recall that by a weak solution to \eqref{eq:Helmholtz II}  we mean that $u\in H^1_{loc}(\mathbb{R}^n)$ and that it satisfies
$$\int_{\mathbb{R}^N}\eta\nabla u\cdot\nabla\varphi-k^2qu\varphi=\int_{\mathbb{R}^n} h\varphi \quad\text{for any }\varphi\in C^{\infty}_0(\mathbb{R}^n).$$
The limit in \eqref{eq:sommerfeldn} has to hold uniformly for every direction $\hat{x}=x/|x|\in \mathbb{S}^{n-1}$ and is also known as the {\it Sommerfeld radiation condition} which characterizes the radiating nature of the scattered wave field $u^s$ (cf. \cite{ColKre,Ned}). There exists a unique weak solution $u(x,d,k)=u^-\chi_{\Omega}+u^+\chi_{\mathbb{R}^n\backslash\overline{\Omega}}\in H_{loc}^1(\mathbb{R}^n)$ to \eqref{eq:Helmholtz II}, and we refer to Appendix in \cite{LSSZ} for a convenient proof. We remark that, if the coefficients are regular enough, say, satisfying those regularity conditions specified at the beginning of the present section, \eqref{eq:Helmholtz II} corresponds to the following transmission problem
\begin{equation}\label{eq:Helmholtz IIclassic}
\begin{cases}
\ \displaystyle{\sum_{i,j=1}^{n}\frac{\partial}{\partial x_i}\left(\eta^{ij}\frac{\partial}{\partial x_j} u^-(x,d,k)\right)+{k^2}q u^-(x,d,k)}=-h\qquad\qquad & x\in\Omega,\\
\ \displaystyle{\Delta u^+(x,d, k)+k^2 u^+(x,d,k)=0}\quad & x\in \mathbb{R}^n\backslash\overline{\Omega},\\
\ \displaystyle{u^-(x)=u^+(x),\ \ \sum_{i,j=1}^n(\nu_i\eta^{ij}\frac{\partial u^-}{\partial x_j})(x)=(\nu\cdot\nabla u^+)(x)}\quad & x\in\partial\Omega,\\
\ \displaystyle{u^+(x,d,k)=e^{i k x\cdot d}+u^s(x,d,k)}\quad & x\in\mathbb{R}^n\backslash\overline{\Omega},\\
\ \displaystyle{\lim_{r\rightarrow\infty}r^{\frac{n-1}{2}}\left(\frac{\partial u^s(x)}{\partial r}-i k u^s(x)\right )=0}\quad & r=|x|,
\end{cases}
\end{equation}
where $\nu=(\nu_i)_{i=1}^n$ is the outward unit normal vector to $\partial\Omega$.

Furthermore, $u(x)$ admits the following asymptotic
development as $|x|\rightarrow+\infty$
\begin{equation}\label{eq:asymptotic}
u(x,d,k)=e^{ik x\cdot d}+\frac{e^{i k |x|}}{|x|^{\frac{n-1}{2}}}\, a_\infty\left(\frac{x}{|x|},d,k\right)+\mathcal{O}\left(\frac{1}{|x|^{\frac{n+1}{2}}}\right).
\end{equation}
In \eqref{eq:asymptotic}, $a_\infty(\hat{x},d,k)$ with $\hat{x}:=x/|x|\in\mathbb{S}^{n-1}$ is known as the {\it far-field pattern} or the {\it scattering
amplitude}, which depends on the impinging direction $d$ and the wave number $k$ of the incident wave $u^i(x):=e^{ik x\cdot d}$, the observation direction $\hat{x}$, and obviously, also the underlying
scattering object $(\Omega;\eta,q, h)$. In the following, we shall also write $a_\infty(\hat{x}, d; (\Omega; \eta, q, h))$ to indicate such dependences, noting that we consider $k$ to be fixed and we drop the dependence on $k$.
An important inverse scattering problem arising in practical applications is to recover the medium $(\Omega; \eta, q)$ or/and the source term $h$ by knowledge of $a_\infty(\hat{x},d)$.  This inverse problem is of fundamental importance in many areas of science and technology, such as radar and sonar, geophysical exploration, non-destructive testing, and medical imaging to name just a few; see \cite{ColKre,Isa} and the references therein. As a general remark, we would like to mention that the acoustic mediums are usually dispersive, and the time-harmonic measurements will provide more accurate reconstruction for the inverse problem. In this context, an acoustic invisibility cloak could be generally introduced as follows.

\begin{defin}\label{def:cloaking device}
Let $\Omega$ and $D$ be bounded domains such that $D\Subset\Omega$.
$\Omega\backslash\overline{D}$ and $D$ represent, respectively, the
cloaking region and the cloaked region. Let $\Gamma$ and $\Gamma'$
be two subsets of $\mathbb{S}^{n-1}$. $(\Omega\backslash\overline{D};
\eta_c, q_c)$ is said to be an (ideal/perfect) {\it invisibility
cloaking device} for the region $D$ if
\begin{equation}\label{eq:definvisibility}
a_\infty\left(\hat{x},d; (\Omega;\eta_e,q_e), h \right)=0\quad \mbox{for}\ \
\hat{x}\in\Gamma,\ d\in\Gamma',
\end{equation}
where the extended object
\[
(\Omega;\eta_e,q_e)=\begin{cases}
\ \eta_a, q_a\quad & \mbox{in\ \ $D$},\\
\ \eta_c,q_c\quad & \mbox{in\ \ $\Omega\backslash\overline{D}$},
\end{cases}
\]
with $(D; \eta_a,q_a)$ denoting a target medium and $h$ denoting a active source/sink inside $D$. If $\Gamma=\Gamma'=\mathbb{S}^{n-1}$, then it is called a {\it
full cloak}, otherwise it is called a {\it partial cloak} with
limited apertures $\Gamma$ of observation angles, and $\Gamma'$ of
impinging angles.
\end{defin}

By Definition~\ref{def:cloaking device}, we have that the cloaking layer $(\Omega\backslash \overline{D}; \eta_c, q_c)$ makes the target medium $(D;\eta_a, q_a)$ together with a source/sink $h$ invisible to the exterior scattering measurements when the detecting waves come from the aperture $\Gamma'$ and the observations are made in the aperture $\Gamma$.

The EM scattering and the corresponding invisibility cloaking can be introduced in a similar manner. Let $\varepsilon=\varepsilon(x)=(\varepsilon^{ij})\in\mathbb{R}^{3\times 3}$, $\mu=\mu(x)=(\mu^{ij}(x))\in\mathbb{R}^{3\times 3}$ and $\sigma=\sigma(x)=(\sigma^{ij}(x))\in\mathbb{R}^{3\times 3}$ be symmetric-matrix valued measurable real functions such that both $\varepsilon$ and $\mu$ are regular (cf. \eqref{eq:reg1}) and $\sigma$ satisfies
\begin{equation}\label{eq:uniform elliptic 2}
0\leq \sum_{i,j=1}^3 \sigma^{ij}(x)\xi_i\xi_j\leq \lambda^{-1} |\xi|^2, \quad \mbox{for any $\xi\in\mathbb{R}^n$ and for a.e. $x\in\mathbb{R}^n$}.
\end{equation}
Physically, functions $\varepsilon$, $\mu$ and $\sigma$ stand respectively for the electric permittivity, magnetic permeability and conductivity tensors of a {\it regular} EM medium. We assume the inhomogeneity of the EM medium is compactly supported such that $\varepsilon=\varepsilon_0$, $\mu=\mu_0$ and $\sigma=\sigma_0$ in $\mathbb{R}^3\backslash\overline{\Omega}$, where $\sigma_0^{ij}=0$ and $\varepsilon_0^{ij}=\mu_0^{ij}=\delta^{ij}$ denote the EM parameters for the homogeneous background space. Let
\begin{equation}\label{eq:EM incident}
E^i(x):=p e^{ikx\cdot d},\quad H^i(x):=\frac{1}{ik}\nabla\wedge E^i(x),\quad x\in\mathbb{R}^3,
\end{equation}
be a pair of time-harmonic EM plane waves. Here, $k\in\mathbb{R}_+$ and $d\in\mathbb{S}^2$ denotes, respectively, the wave number and the impinging direction, and $p\in\mathbb{R}^3$ with $p\perp d$ denotes the polarization of the plane waves. 
Then the EM wave propagation in the whole space $\mathbb{R}^3$ with 
an EM medium inclusion $(\Omega; \varepsilon,\mu,\sigma)$ 
as described above 
is governed by the following Maxwell system
\begin{equation}\label{eq:Maxwell whole}
\begin{cases}
\displaystyle{\nabla\wedge E-i\omega \mu H=0}\qquad\qquad &\mbox{in\ \ $\mathbb{R}^3$,}\\
\displaystyle{\nabla\wedge H+i\omega\left(\varepsilon+i\frac{\sigma}{\omega}\right) E={J}}\quad &\mbox{in\ \ $\mathbb{R}^3$},\\
\displaystyle{ E^-= E|_{\Omega},\quad {E}^+=({E}-E^i)|_{\mathbb{R}^3\backslash\overline{\Omega}}},\\
\displaystyle{ H^-= H|_{\Omega},\quad {H}^+=({H}-H^i)|_{\mathbb{R}^3\backslash\overline{\Omega}}},\\
\displaystyle{\lim_{|x|\rightarrow+\infty}|x|\left| (\nabla\wedge E^+)(x)\wedge\frac{x}{|x|}-i\omega E^+(x) \right|=0},
\end{cases}
 \end{equation}
where ${J}\in\mathbb{C}^3$ denotes an electric current density supported in 
$\Omega$ and $J\in L^2(\Omega)^{3}$. In \eqref{eq:Maxwell whole}, $E$ and $H$ are respectively the electric and magnetic fields, and $E^+$ and $H^+$  are the scattered fields. The last relation 
in \eqref{eq:Maxwell whole} is called the Silver-M\"uller radiation condition, which characterizes the radiating nature of the scattered wave fields $E^+$ and $H^+$. For a regular EM medium $(\Omega; \varepsilon, \mu, \sigma)$ and an active electric current 
$J\in L^2(\Omega)^{3}$, there exists a unique pair of solutions $E, H\in 
H_{loc}(\nabla\wedge; \mathbb{R}^3)$ (see \cite{Lei,Ned}).
Here and in what follows 
\[
H_{loc}(\nabla\wedge; X)=\{ U|_B\in H(\nabla\wedge; B)|\ B\ \ \mbox{is any bounded subdomain of $X$} \}
\]
and
\[
H(\nabla\wedge; B)=\{ U\in (L^2(B))^3|\ \nabla\wedge U\in (L^2(B))^3 \}.
\]
Furthermore, $E^+$ admits the asymptotic expression as $|x|\rightarrow \infty$ (cf. \cite{ColKre}):
\begin{equation}\label{eq:farfield EM}
E^+(x)=\frac{e^{i\omega |x|}}{|x|} A_\infty\left(\frac{x}{|x|}; p,d,k\right)+\mathcal{O}\left(\frac{1}{|x|^2}\right)
\end{equation}
where $A_\infty(\hat{x}; p,d,k)$ is known as the EM {\it scattering amplitude}. 
In the sequel, we shall also write $A_\infty(\hat{x}; E^i, (\Omega; \varepsilon, \mu, \sigma), J)$ to specify the related dependences.
The inverse EM scattering problem is to recover $(\Omega;\varepsilon,\mu,\sigma)$ and/or $J$ by knowledge of $A_\infty(\hat{x},d)$. The partial- and full-cloaks in the EM scattering could be introduced in a completely similar manner to Definition~\ref{def:cloaking device}. 

\section{Transformation acoustics and electromagnetics}

Let $\Omega$ and $\widetilde\Omega$ be two bounded Lipschitz domains, and suppose there exists a bi-Lipschitz and orientation-preserving mapping 
$$
\widetilde x=F(x): \Omega\rightarrow\widetilde\Omega.
$$ 
The key ingredients of the transformation acoustics are summarized in the following lemma.

\begin{lem}\label{lem:trans acoustics}

Let $(\Omega; \eta, q, h)$ be a scattering configuration supported in $\Omega$. Define the {\it push-forwarded} scattering configuration as follows,
\begin{equation}\label{eq:pushforward}
(\widetilde\Omega; \widetilde\eta, \widetilde q, \widetilde h):=F_*(\Omega; \eta, q, h),
\end{equation}
where
\begin{equation}
\begin{split}
&\widetilde{\eta}(\widetilde
x)=F_*\eta(\widetilde x):=\frac{DF(x)\cdot \eta(x)\cdot DF(x)^T}{|\mbox{\emph{det}}(DF(x))|}\bigg|_{x=F^{-1}(\widetilde x)},\\
&\widetilde{q}(\widetilde x)=F_*q(\widetilde x):=\frac{q(x)}{|\mbox{\emph{det}}(DF(x))|}\bigg|_{x=F^{-1}(\widetilde x)},
\end{split}
\end{equation}
and
\begin{equation}
\begin{split}
\widetilde{h}(\widetilde{x})=F_*h(\widetilde x):=\frac{h(x)}{|\mbox{\emph{det}}(DF(x))|}\bigg|_{x=F^{-1}(\widetilde x)},
\end{split}
\end{equation}
where $DF$ denotes the Jacobian matrix of $F$. 
Then $u\in H^1(\Omega)$ solves the Helmholtz equation
\[
\nabla\cdot(\eta\nabla u)+k^2 q u=-h\quad \mbox{in\ $\Omega$},
\]
if and only if the pull-back field $\widetilde u=(F^{-1})^*u:=u\circ F^{-1}\in H^1(\widetilde\Omega)$ solves
\[
\widetilde\nabla\cdot(\widetilde\eta\nabla \widetilde u)+k^2\widetilde q \widetilde u=-\widetilde h\quad \mbox{in\ $\widetilde\Omega$}.
\]
We have made use of $\nabla$ and $\widetilde\nabla$ to distinguish the
differentiations respectively in $x$- and $\widetilde x$-coordinates.
\end{lem}

A convenient proof of Lemma~\ref{lem:trans acoustics} can be found, e.g., in \cite{KOVW}. As a direct consequence of Lemma~\ref{lem:trans acoustics}, one can show that if $F:\overline\Omega\rightarrow\overline\Omega$ is a bi-Lipschitz mapping such that $F|_{\partial\Omega}=\text{Id}$, then
\begin{equation}\label{eq:trans acoustics equal}
a_\infty(\hat{x}; (\Omega;\eta, q, h))=a_\infty(\hat{x}; F_*(\Omega; \eta, q, h)).
\end{equation}

The key ingredients of the transformation electromagnetics are summarized in the following lemma.

\begin{lem}\label{lem:trans opt}
Let $(\Omega;\varepsilon,\mu,\sigma, J)$ be an EM scattering configuration supported in $\Omega$. Define the {\it push-forwarded} scattering configuration as follows,
\begin{equation}\label{eq:pushforward EM}
(\widetilde\Omega; \widetilde\varepsilon,\widetilde\mu,\widetilde\sigma, \widetilde J):=F_*(\Omega; \varepsilon,\mu,\sigma,J),
\end{equation}
where
\begin{equation}
\begin{split}
&\widetilde{m}(\widetilde
x)=F_*m(\widetilde x):=\frac{DF(x)\cdot m(x)\cdot DF(x)^T}{|\mbox{\emph{det}}(DF(x))|}\bigg|_{x=F^{-1}(\widetilde x)},\quad m=\varepsilon,\mu\ \mbox{and}\ \sigma,\\
&\widetilde J(\widetilde x)=F_* J(\widetilde x):=\frac{1}{|\text{\emph{det}}(DF(x))|}(DF(x)) J(x)\bigg|_{x=F^{-1}(\widetilde x)},
\end{split}
\end{equation}
Suppose that $E, H\in H(\nabla\wedge;\Omega)$ are the EM fields satisfying
\[
\begin{split}
\nabla\wedge E-ik \mu H=0\qquad & \mbox{in\ \ $\Omega$,}\\
\nabla\wedge H+ik \left( \varepsilon+i\frac{\sigma}{\omega} \right) E=J\qquad & \mbox{in\ \ $\Omega$},
\end{split}
\]
If we define the {\emph pull-back fields} by
\[
\begin{split}
\widetilde E &=({F}^{-1})^* E:=(D{F})^{-T}E\circ {F}^{-1},\\
\widetilde H &=({F}^{-1})^* H:=(D {F})^{-T}H\circ {F}^{-1},
\end{split}
\]
then the pull-back fields $\widetilde E, \widetilde H\in H(\widetilde\nabla\wedge; 
\widetilde\Omega)$ 
satisfy the following Maxwell equations
\[
\begin{split}
\widetilde\nabla \wedge \widetilde E-ik \widetilde\mu \widetilde H=& 0\qquad \mbox{in\ \ $\widetilde\Omega$},\\
\widetilde\nabla \wedge \widetilde H+ik \left(\widetilde\varepsilon+i\frac{\widetilde\sigma}{k}\right) \widetilde E=& \widetilde J\quad\ \ \mbox{in\ \ $\widetilde\Omega$}.
\end{split}
\]
\end{lem}

The proof of Lemma~\ref{lem:trans opt} is given in \cite{LiuZhou}. As a direct consequence of Lemma~\ref{lem:trans opt}, one can show that if $F:\overline\Omega\rightarrow\overline\Omega$ is a bi-Lipschitz mapping such that $F|_{\partial\Omega}=\text{Id}$, then
\begin{equation}\label{eq:trans em equal}
A_\infty(\hat{x}; (\Omega;\varepsilon,\mu,\sigma,J))=A_\infty(\hat{x}; F_*(\Omega; \varepsilon,\mu,\sigma,J)).
\end{equation}

\section{Regularized cloaks in acoustic scattering}

In this section, we discuss the results on regularized cloaks in acoustic scattering governed by the Helmholtz equation. 

\subsection{Regularized full cloaks}\label{sect:41}

In the rest of this paper, we let $\Omega$ and $D$ be two bounded simply connected smooth domains in $\mathbb{R}^N$ containing the origin such that $D\Subset\Omega$. Define
\begin{equation}\label{eq:scale domain}
D_\rho:=\{\rho x; x\in D\},\quad \rho\in\mathbb{R}_+. 
\end{equation}
Throughout, we assume there exists a bi-Lipschitz and orientation-preserving mapping,
\begin{equation}\label{eq:F1}
F_\rho^{(1)}: \overline{\Omega}\backslash D_\rho\rightarrow\overline\Omega\backslash D,\quad F_\rho^{(1)}|_{\partial\Omega}=\mbox{Identity},
\end{equation}
for $0<\rho<1$. That is, $F_\rho^{(1)}$ blows up $D_\rho$ within $\Omega$ (see Fig.~5). Let
\begin{equation}\label{eq:F2}
F_\rho^{(2)}(x)=\frac{x}{\rho},\quad x\in D_\rho,
\end{equation}
and
\begin{equation}\label{eq:Fwhole}
F_\rho=\begin{cases}
F_\rho^{(1)}\qquad & \mbox{on}\ \ \overline\Omega\backslash D_\rho,\\
F_\rho^{(2)}\qquad & \mbox{on}\ \ D_\rho. 
\end{cases}
\end{equation}
Consider a virtual scattering configuration as follows
\begin{equation}\label{eq:virtual a1}
\Omega; \eta, q=\begin{cases}
I, 1\qquad & \mbox{on}\ \ \Omega\backslash D_\rho,\\
\eta_l, q_l\qquad & \mbox{on}\ \ D_\rho\backslash D_{\rho/2},\\
\eta_a, q_a\qquad & \mbox{on}\ \ D_{\rho/2},
\end{cases}
\end{equation}
where $(\eta_l, q_l)$ and $(\eta_a, q_a)$ shall be specified in the sequel. Let $(\Omega;\widetilde\eta, \widetilde q)$ be a physical scattering configuration given by
\begin{equation}\label{eq:physical a1}
(\Omega; \widetilde\eta, \widetilde q):=(F_\rho)_*(\Omega; \eta, q),
\end{equation}
where $F_\rho$ is given in \eqref{eq:Fwhole} and $(\Omega; \eta, q)$ is given in \eqref{eq:virtual a1}. 

\begin{figure}[htbp]\label{fig:5}
\begin{center}
\includegraphics[width=.7\linewidth]{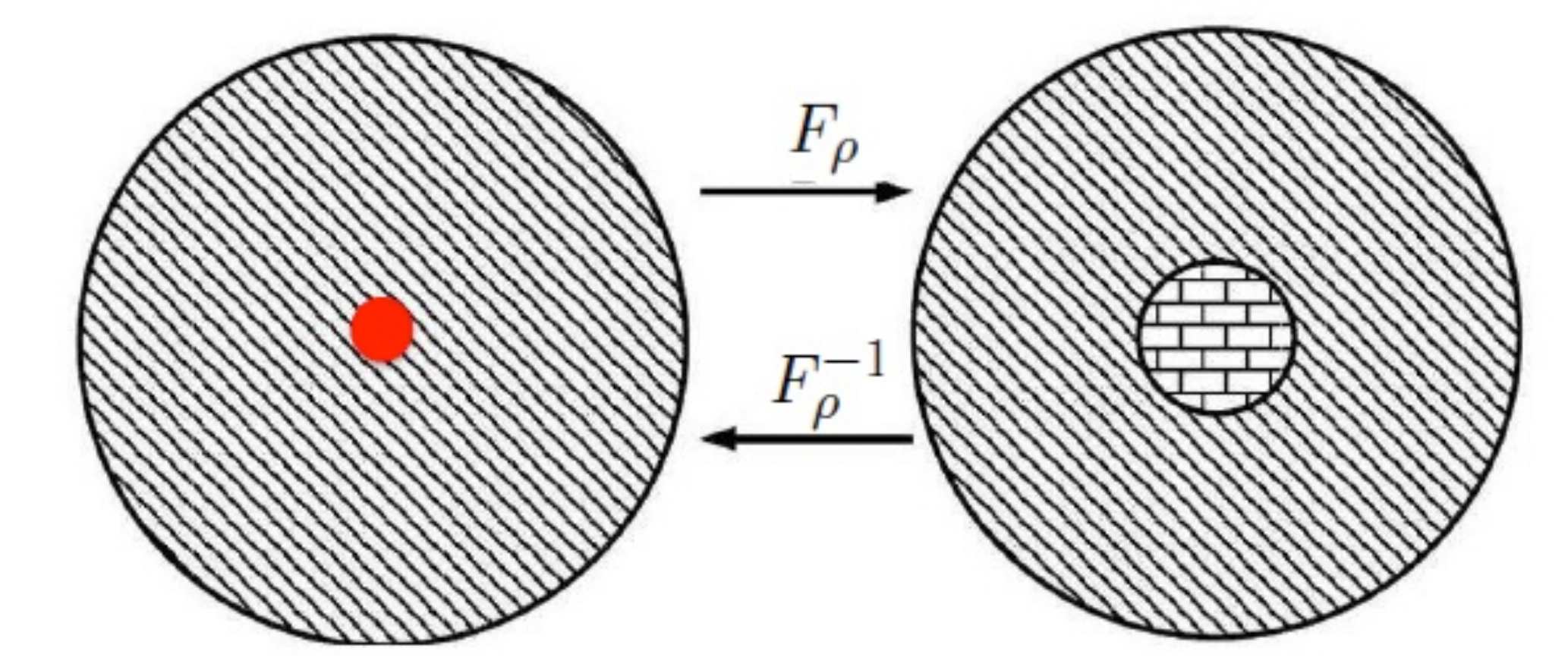}
\caption{Map $F_\rho^{(1)}:B(0,2)\setminus \overline B(0,\rho)\to B(0,2)\setminus\overline
B(0,1)$}\end{center}\end{figure}

\begin{thm}[\cite{KOVW} and \cite{LiuSun}]\label{thm:cloakf1}
Let $(\Omega;\widetilde\eta, \widetilde q)$ be a scattering configuration given in \eqref{eq:physical a1}. There exists $\rho_0\in\mathbb{R}_+$ such that when $\rho<\rho_0$:
\begin{enumerate}[i).]
\item If
\begin{equation}\label{eq:highloss}
\eta_l(x)=I,\quad q_l(x)=1+i\rho^{-2},\quad x\in D_{\rho}\backslash D_{\rho/2},
\end{equation}
then
\begin{equation}\label{eq:highloss2}
\|a_\infty(\hat{x}; (\Omega; \widetilde\eta, \widetilde q))\|_{L^\infty(\mathbb{S}^{n-1})}\leq C\|u^i\|_{L^2(B_R)}\times\begin{cases} |\ln\rho|^{-1}, & \quad n=2,\\ \rho, &\quad n=3. \end{cases}
\end{equation}

\item If
\begin{equation}\label{eq:highdens}
\eta_l(x)=\rho^2I,\quad q_l(x)=1+i,\quad x\in D_{\rho}\backslash D_{\rho/2},
\end{equation}
then
\begin{equation}\label{eq:highdens2}
\|a_\infty(\hat{x}; (\Omega; \widetilde\eta, \widetilde q))\|_{L^\infty(\mathbb{S}^{n-1})}\leq C\rho^n\|u^i\|_{L^2(B_R)}.
\end{equation}
\end{enumerate}
In both the estimates \eqref{eq:highloss2} and \eqref{eq:highdens2}, $B_R$ denotes a central ball containing $\Omega$, and the generic constant $C$ is independent of $\rho$, $u^i$, $\eta_a$ and $q_a$. 
\end{thm}

In the sequel, we let
\begin{equation}\label{eq:cloak1}
(\Omega\backslash D; \widetilde\eta_c^\rho, \widetilde q_c^\rho):=(F_\rho)_*(\Omega\backslash D_\rho; I, 1),
\end{equation}
and
\begin{equation}\label{eq:cloak2}
(D\backslash D_{1/2}; \widetilde\eta_l,\widetilde q_l):=(F_\rho)_*(D_\rho\backslash D_{\rho/2}; \eta_l, q_l).
\end{equation}
In the theoretic limit case $\rho=+0$, $D_\rho$ degenerate to a singular point, and $F_\rho$ blows up the singular point to $D$ within $\Omega$. Moreover, $(\Omega\backslash D; \widetilde\eta_c^0, \widetilde q_c^0)$ is the ideal cloaking layer considered in \cite{GKLU3} which can be used to cloak an arbitrary (but regular) passive medium $(D; \widetilde\eta_a, \widetilde q_a)$. Theorem~\ref{thm:cloakf1} indicates that the regularized cloaking layer $(\Omega\backslash D; \widetilde \eta_c^\rho, \widetilde q_c^\rho)$ together with the lossy layer $(D\backslash D_{1/2}; \widetilde\eta_l, \widetilde q_l)$ produces an approximate cloaking device within accuracy $e(\rho)$, which can be used to cloak an arbitrary passive medium $(D_{1/2}; \widetilde\eta_a,\widetilde q_a)=(F_\rho)_*(\eta_a, q_a)$. Here, $e(\rho)$ denotes the RHS $\rho$-terms in \eqref{eq:highloss2} and \eqref{eq:highdens2}. It is remarked that the lossy layer $(D\backslash D_{1/2}; \widetilde\eta_l, \widetilde q_l)$ is necessary for achieving a practical near-cloaking device since otherwise there exist cloak-busting inclusions $(D;\eta_\rho, q_\rho)$ which defy any attempt to nearly cloak them (see \cite{KOVW}). In \cite{LiLiuSun}, \eqref{eq:highloss} is referred to as a {\it high-loss layer}, and \eqref{eq:highdens} is referred to as a {\it high-density layer}. The high-lossy layer is shown to be a finite realization of a sound-soft layer, and the high-density layer is shown to be finite realization of a sound-soft layer. Both the estimates in Theorem~\ref{thm:cloakf1} are shown to be sharp for the respective constructions. Here, we would like to mention a few words about the arguments in deriving the cloaking assessments. By using Lemma~\ref{lem:trans acoustics}, one has 
\begin{equation}\label{eq:tt1}
a_\infty(\hat{x}; (\Omega; \widetilde\eta, \widetilde q))=a_\infty(\hat{x}; (\Omega; \eta, q)),
\end{equation}
where $(\Omega; \eta, q)$ is given in \eqref{eq:virtual a1}. The scattering configuration given in \eqref{eq:virtual a1} is actually a small inclusion supported in $D_\rho$ with an arbitrary content in $D_{\rho/2}$ enclosed by a thin lossy layer $(D_\rho\backslash D_{\rho/2}; \eta_l, q_l)$. Hence, in order to assess the corresponding cloaking constructions, it suffices for one to consider the scattering estimate due to a small inclusion possessing the peculiar structure as described above. 

In \cite{Liu2}, the cloaking of active contents by employing a general lossy layer is considered. We let
\begin{equation}\label{eq:gl1}
\eta_l(x)=\rho^{r}\gamma(x/\rho),\quad x\in D_\rho\backslash D_{\rho/2},
\end{equation}
where $r>2-N/2$ and $\gamma(x)\in C^2(\overline D\backslash D_{1/2})$ is positive function that is bounded below, and 
\begin{equation}\label{eq:gl2}
q_l(x)=\alpha(x/\rho)+i\beta(x/\rho),\quad x\in D_\rho\backslash D_{\rho/2},
\end{equation}
where $\alpha$ and $\beta$ are positive functions that are bounded below and above. 

\begin{thm}[\cite{Liu2}]\label{thm:cloakf2}
Let $(\Omega; \eta, q)$ be a virtual scattering configuration as that given in \eqref{eq:virtual a1} but with $(D_\rho\backslash D_{\rho/2}; \eta_l,q_l)$ given in \eqref{eq:gl1}--\eqref{eq:gl2}, and let $(\Omega; \widetilde\eta,\widetilde q)$ be the corresponding physical scattering configuration given in \eqref{eq:physical a1}. Let $h\in L^2(D_{1/2})$ be a source/sink in the physical scattering configuration. Assume that
\begin{equation}\label{eq:assum1}
\Im \widetilde q_a\geq \gamma_0>0\quad \mbox{on}\ \ supp(h)\subset D_{1/2}.
\end{equation}
Then there exists $\rho_0\in\mathbb{R}_+$ such that when $\rho<\rho_0$
\begin{equation}\label{eq:active est}
\begin{split}
\|a_\infty(\hat{x};(\Omega; \widetilde\eta,\widetilde q))\|_{L^\infty(\mathbb{S}^{n-1})} \leq C & \bigg(\rho^{\min(n+2r-4,n)}\|u^i\|_{L^2(B_R)}\\
&+\rho^{\min(n/2,n/2+r-2)}\|h\|_{L^2(D_{1/2})} \bigg),
\end{split}
\end{equation}
where $C$ is independent of $\rho$, $u^i$, $\eta_a$, $\Re q_a$, $r$ and $h$. 
\end{thm}

The general lossy layer considered in \eqref{eq:gl1}--\eqref{eq:gl2} will be of interest when production fluctuation occurs. Theorem~\ref{thm:cloakf2} indicates the effective way of cloaking active contents. In particular, one should maintain an absorbing environment for the place where the the source/sink is located, and this viewpoint shall also be adopted in our subsequent study on the partial cloaks in acoustic scattering. Finally, we give a brief discussion on the idea of proving Theorem~\ref{thm:cloakf2}. As discussed earlier, by virtue of \eqref{eq:tt1}, one suffices to estimates the scattering due to a small inclusion supported in $D_\rho$ with an arbitrary passive content and an active source in $D_{\rho/2}$ enclosed by a thin lossy layer $(D_\rho\backslash D_{\rho/2}; \eta_l, q_l)$ in \eqref{eq:gl1}-\eqref{eq:gl2}. First, by using the structure of the lossy layer and a variational argument, one can control the energy of the wave field in the thin layer $D_\rho\backslash D_{\rho/2}$. Then, by a duality argument, one can control the trace of the wave velocity on $\partial D_{\rho}$. Hence, the problem is reduced to estimating the scattering from a small inclusion $D_\rho$ with a prescribed trace on $\partial D_\rho$, and this can be achieved by the layer-potential technique.

\subsection{Regularized partial cloaks}\label{sect:42}

In \cite{LiLiuRonUhl}, the regularized blow-up construction for acoustic cloaks has been extended to an extremely general setting, which we shall describe in this section. As before, we start our discussion in the virtual space. In the sequel, we describe a typical but a bit simplified virtual scattering configuration in \cite{LiLiuRonUhl}.

Let $K_0$ be a compact subset in $\mathbb{R}^N$ such that $G:=\mathbb{R}^n\backslash K_0$ is connected.
We denote by  $d:\mathbb{R}^n\to[0,+\infty)$ the distance function from $K_0$ defined as follows
$$d(x)=\mathrm{dist}(x,K_0)\quad\text{for any }x\in\mathbb{R}^n.$$ Assume that there exists a Lipschitz function $\widetilde{d}:\mathbb{R}^n\to [0,+\infty)$ such that
the following properties are satisfied.

First, there exist constants $a$ and $b$, $0<a\leq1\leq b$, such that
$$ad(x)\leq\widetilde{d}(x)\leq bd(x)\quad\text{for any }x\in\mathbb{R}^n.$$

For any $\rho>0$, let $K_{\rho}=\{x\in\mathbb{R}^n:\ \widetilde{d}(x)\leq \rho\}$.
For some constants $\rho_0>0$, $p>2$, $C_1>0$ and $R>0$, we require that for any $\rho\in\mathbb{R}_+$, $\rho\leq \rho_0$,  $K_{\rho}\subset \overline{B_R}$,
$\mathbb{R}^n\backslash K_{\rho}$ is connected
and
$$\|u\|_{L^p(B_{R+1}\backslash K_{\rho})}\leq C_1\|u\|_{H^1(B_{R+1}\backslash K_{\rho})}\quad\text{for any }u\in H^1(B_{R+1}\backslash K_{\rho}).$$

We notice that, a simple sufficient condition for these assumptions to hold is that $K_0$ is a compact convex set. Another sufficient condition is that $K_0$ is a Lipschitz scatterer. Roughly speaking, $K_0$ is a Lipschitz scatterer if $\partial K_0$ is consisting of Lipschitz hypersurfaces, and we refer to \cite[Section~4]{Men-Ron} for detailed discussion. Finally, $K_0$ may be the union of a finite number of pairwise disjoint
compact convex sets and Lipschitz scatterers. 

Now, we consider a virtual scattering configuration $(\mathbb{R}^n; \eta^\rho, q^\rho, h^\rho)$ satisfying the following assumptions. 

\begin{enumerate}[a)]
\item We have that $\eta^{\rho}(x)=I$ and $q^{\rho}(x)=1$ for almost every $x\in \mathbb{R}^n\backslash K_{\rho}$.
\item There exist a continuous nondecreasing function $\omega_1:(0,\rho_0]\to (0,+\infty)$, such that
$\lim_{s\to 0^+}\omega_1(s)=0$, and positive constants $E_1$, $E'_1$, $\Lambda$, and $E_2$ such that for almost any $x\in K_{\rho}\backslash K_{\rho/2}$
\begin{equation}\label{eq:ll1}
0<\Re q^{\rho}(x)\leq E_1(\omega_1(\rho))^{-1}\quad\text{and}\quad 0<E'_1(\omega_1(\rho))^{-1}\leq \Im q^{\rho}(x).
\end{equation}
Furthermore
\begin{equation}\label{firstcondsigma}
\eta^{\rho}(x)\xi\cdot\xi\leq \Lambda\|\xi\|^2\quad\text{for any }\xi\in\mathbb{R}^n\text{ and for a.e. }x\in K_{\rho}\backslash K_{\rho/2}
\end{equation}
and
\begin{equation}\label{secondcondsigmabis}
\frac{1}{\rho^2}\eta^{\rho}(x)\nabla \widetilde{d}(x)\cdot\nabla \widetilde{d}(x)\leq E_2\quad\text{for a.e. }x\in K_{\rho}\backslash K_{\rho/2}.
\end{equation}
Finally we require that
\begin{equation}\label{secondassumption}
\lim_{\rho\to 0^+}\int_{K_{\rho}\backslash K_{\rho/2}}|q^{\rho}|=0.
\end{equation}

\item 
$h^{\rho}(x)=0$ for almost every $x\in \mathbb{R}^N\backslash K_{\rho/2}$. Moreover, there exists a positive constant $E_3$ such that
\begin{equation}\label{sourceassumptionbis}
\int_{K_{\rho/2}}(k^2\Im q^{\rho})^{-1}|h^{\rho}|^2 \leq E_3.
\end{equation}
Notice that the above condition means that $\Im q^{\rho}>0$ almost everywhere on the set $\{h^{\rho}\neq 0\}$ and that the integral is actually performed on such a set.

\end{enumerate}

For the condition \eqref{sourceassumptionbis}, we note that it remains unchanged under push-forward due to the fact that
\[
\int_{\Omega} q^{-1} |h|^2=\int_{\widetilde\Omega} \widetilde q^{-1}|\widetilde h|^2,
\] 
where $(\widetilde\Omega; \widetilde q, \widetilde h)=F_*(\Omega; q, h)$.

Let $u_\rho$ denote the scattering solution corresponding to the scattering configuration $(\mathbb{R}^N; \eta^\rho, q^\rho, h^\rho)$ as described above; see \eqref{eq:Helmholtz II} and \eqref{eq:Helmholtz IIclassic}. The following significant convergence result is proved in \cite{LiLiuRonUhl}. 

\begin{thm}[\cite{LiLiuRonUhl}]\label{scatteo}
Under the previous assumptions, $v_\rho=u_\rho(1-\chi_{K_{\rho/2}})$ converges to a function $u$ strongly in $L^2(B_r)$ for any $r>0$, with $u$ solving
\begin{equation}\label{uscateq}
\left\{\begin{array}{ll}
\Delta u+k^2u=0 &\text{in }\ \mathbb{R}^n\backslash K_0,\\
u=u^i+u^s &\text{in }\ \mathbb{R}^n\backslash K_0,\\
\nabla u\cdot \nu=0 &\text{on }\ \partial K_0,\\
\lim_{r\to+\infty}r^{(n-1)/2}\left(\frac{\partial u^s}{\partial r}-i ku^s\right)=0 &r=|x|,
\end{array}\right.
\end{equation}
and $u=0$ in $K_0$.
\end{thm}

Based on Theorem~\ref{scatteo}, one can construct regularized full- or partial-cloaks through blow-up construction via pushing forward the virtual scattering configuration $(\mathbb{R}^n; \eta^\rho, q^\rho)$. We assume that there exists a bi-Lipschitz and orientation-preserving mapping $F_\rho$ such that
\begin{equation}\label{eq:ff1}
F_\rho=\mbox{Identity\ \ on}\ \ \mathbb{R}^n\backslash\overline{\Omega};\quad F_\rho^{(1)}\ \ \mbox{on}\ \ \overline{\Omega}\backslash K_\rho;\quad F_\rho^{(2)}\ \ \mbox{on}\ \ K_\rho
\end{equation}
with
\begin{equation}\label{eq:ff2}
F_\rho^{(2)}(K_{\rho/2})=\overline{D}_{1/2},\quad F_\rho^{(2)}(K_\rho)=\overline{D},\quad F_\rho^{(1)}(\Omega\backslash K_\rho)=\overline{\Omega}\backslash\overline{D}. 
\end{equation}

We first consider the construction of regularized full-cloaks. Set
\begin{equation}\label{eq:ppvv}
 K_\rho=D_{\rho},\ \ K_{\rho/2}=D_{\rho/2},
\end{equation}
and
\begin{align}
& \widetilde\eta_c^\rho=(F_\rho^{(1)})_*I,\quad \widetilde q_c^\rho=(F_\rho^{(1)})_* 1\quad\mbox{in\ $\Omega\backslash D$},\label{eq:cl1}\\
&(D\backslash {D}_{1/2};\widetilde\eta^\rho_l, \widetilde q_l^\rho)=(F_\rho)_*(K_{\rho}\backslash {K}_{\rho/2}; \eta_l^\rho, q_l^\rho),\label{eq:ly1}
\end{align}
with
\begin{equation}\label{eq:ly2}
\eta_l^\rho=c_1(x)\rho^{2},\quad q_l^\rho=(c_2(x)+i c_3(x))\rho^{-n+1}.
\end{equation}
In \eqref{eq:ly2}, $c_1(x)$ is a symmetric-matrix valued measurable function, and $c_2(x), c_3(x)$ are bounded real valued measurable function such that
\begin{equation}\label{eq:ll2}
\lambda_0|\xi|^2\leq c_1(x)\xi\cdot\xi\leq \Lambda_0|\xi|^2,\quad \lambda_0\leq c_2(x), c_3(x)\leq \Lambda_0\qquad\mbox{for a.e. $x\in K_\rho\backslash K_{\rho/2}$},
\end{equation}
where $\lambda_0$ and $\Lambda_0$ are two positive constants independent of $\rho$. $(D\backslash D_{1/2}; \widetilde\eta_l^\rho, \widetilde q_l^\rho)$ is the chosen lossy layer for the cloaking scheme. By Lemma~\ref{lem:trans acoustics}, it is straightforward to verify that
\begin{equation}\label{eq:ll3}
\widetilde\eta_l^\rho(x)=c_1(\rho x)\rho^n,\quad \widetilde q_l^\rho(x)=(c_2(\rho x)+ic_3(\rho x))\rho,\quad x\in D\backslash D_{1/2}.
\end{equation}

In summarizing the above description, we have a physical scattering configuration as follows
\begin{equation}\label{eq:rc2}
(\mathbb{R}^n; \widetilde{\eta}^\rho, \widetilde{q}^\rho)=\begin{cases}
I, 1\qquad &\ \mbox{in\quad $\mathbb{R}^n\backslash\Omega$},\\
\widetilde{\eta}_c^\rho, \widetilde{q}_c^\rho\qquad & \ \mbox{in\quad $\Omega\backslash D$},\\
\widetilde\eta_l^\rho, \widetilde q_l^\rho\qquad & \ \mbox{in\quad $D\backslash D_{1/2}$},\\
\widetilde\eta_a, \widetilde q_a\qquad &\ \mbox{in\quad $D_{1/2}$},
\end{cases}
\end{equation}

\begin{prop}[\cite{LiLiuRonUhl}]\label{prop:full1}
There exist $\rho_0>0$ and a function $\omega:(0,\rho_0]\to (0,+\infty]$ with $\lim_{s\to 0^+}\omega(s)=0 $, which is independent of $\widetilde\eta_a, \widetilde q_a$ such that for any $\rho<\rho_0$
\begin{equation}\label{eq:full1}
\|a_\infty(\hat{x}; (\mathbb{R}^n; \widetilde\eta^\rho, \widetilde q^\rho))\|_{L^\infty(\mathbb{S}^{n-1})}\leq \omega(\rho).
\end{equation}
\end{prop}

Proposition~\ref{prop:full1} indicates that one has an approximate full invisibility cloak for the construction \eqref{eq:ppvv}--\eqref{eq:rc2}. By the transformation acoustics, one sees that the scattering corresponding to the physical configuration is the same as that corresponding to the virtual configuration, namely $a_\infty(\hat{x}; (\mathbb{R}^n; \widetilde\eta^\rho, \widetilde q^\rho))=a_\infty(\hat{x}; (\mathbb{R}^n; \eta^\rho, q^\rho))$. In the virtual space, the scattering object is supported in $K_\rho$, which degenerates to a single point as $\rho\rightarrow +0$.
The essential point in Proposition~\ref{prop:full1} is that a single point has zero capacity. We recall that $H^1(D)=H^1(D\backslash K_0)$ for any open set $D$ and any compact $K_0\subset D$ with zero capacity.
We emphasize that by following the same spirit, and using Theorem~\ref{scatteo}, one could have more approximate full invisibility cloaks. For example, in $\mathbb{R}^3$, a line segment is also of zero capacity, and hence one can achieve an approximate full cloak by blowing up a `line-segment-like' region in $\mathbb{R}^3$, namely $K_0$ is a line segment; or by blowing up a finite collection of `point-like' and `line-segment-like' regions. We refer to \cite{LiLiuRonUhl} for more concrete constructions. 

Next, we consider the cloaking of active contents. In a completely similar manner, one can show that

\begin{prop}[\cite{LiLiuRonUhl}]\label{prop:full2}
Under the same assumptions of Proposition~\textnormal{\ref{prop:full1}}, we further assume there exists a physical source/sink term
$h\in L^2(D_{1/2})$. Moreover, we require that
\begin{equation}\label{eq:abs cond}
\Im \widetilde{q}_a\geq \lambda_0>0\quad\mbox{on\ $supp(h)$},
\end{equation}
where $\lambda_0$ is a constant.

Then there exist $\varepsilon_0>0$ and a function $\omega:(0,\varepsilon_0]\to (0,+\infty]$ with $\lim_{s\to 0^+}\omega(s)=0 $, which is independent of $\widetilde\sigma_a, \Re\widetilde q_a, \Im\widetilde q_a|_{\{h=0\}}$, such that for any $\rho<\rho_0$
\begin{equation}\label{eq:full2}
\|a_\infty(\hat{x}; (\mathbb{R}^n; \widetilde\eta^\rho, \widetilde q^\rho, h))\|_{L^\infty(\mathbb{S}^{n-1})}\leq \omega(\rho).
\end{equation}

\end{prop}

Compared to the regularized full cloaks discussed in Section~\ref{sect:41}, we note the following two facts. First, the lossy layer employed in Propositions~\ref{prop:full1} and \ref{prop:full2} are much more general than those discussed in Section~\ref{sect:41}. Nonetheless, the general lossy layers proposed in \cite{LiLiuRonUhl} could not include those lossy layers proposed in \cite{KOVW,LiuSun} as special cases. Second, Theorem~\ref{scatteo} only gives the convergence of the scattered wave fields corresponding to the cloaking constructions, and it does not provide the corresponding estimates of the rate of convergence. The rate of convergence would indicate the degree of approximation of the near-cloak to the ideal cloak. Those are interesting issues for further investigation. 

Next, we consider the regularized partial cloaks. The construction of partial cloaking devices will rely on blowing up `partially' small regions in the virtual space.
Let
\begin{equation}\label{eq:k0 2d}
K_0:=\{-a\leq x_1\leq a\}\times \{x_2=0\}\quad\mbox{in\ \ $\mathbb{R}^2$},
\end{equation}
and
\begin{equation}\label{eq:k0 3d}
K_0:=\{-a\leq x_1\leq a\}\times\{ -b\leq x_2\leq b \}\times\{x_3=0\}\quad\mbox{in\ \ $\mathbb{R}^3$}.
\end{equation}

It is noted that $\nu=(0,1)$ in 2D and $\nu=(0,0,1)$ in 3D for $K_0$.
Let $0\leq\tau\leq 1$ and define
\begin{equation}\label{eq:unit sphere set2}
\mathcal{N}_\tau:=\{\theta\in\mathbb{S}^{N-1}:\
|\nu\cdot \theta|\leq\tau\}.
\end{equation}


Next, we consider the scattering problem \eqref{uscateq} with $K_0$ given in \eqref{eq:k0 2d} and \eqref{eq:k0 3d}, which is known as the {\it screen problem} in the literature.

\begin{prop}[\cite{LiLiuRonUhl}]\label{prop:partial virtual}
Let $K_0$ be given in \eqref{eq:k0 2d} and \eqref{eq:k0 3d}, and $\mathcal{N}_\tau$ be given in 
\eqref{eq:unit sphere set2}. Let $u\in H_{loc}^1(\mathbb{R}^n\backslash K_0)$ be the solution to \eqref{uscateq} with $u^i(x)=e^{ikx\cdot d}$. Then there exists a constant $C$, depending only on $a, b$ and $k$, such that
\begin{equation}\label{eq:partial virtual 1}
|a_\infty(\hat{x}, d)|\leq C \tau\quad \mbox{for $\hat{x}\in\mathbb{S}^{n-1}$ and $d\in \mathcal{N}_\tau$},
\end{equation}
and
\begin{equation}\label{eq:partial virtual 2}
|a_\infty(\hat{x}, d)|\leq C \tau\quad \mbox{for $\hat{x}\in\mathcal{N}_\tau$ and $d\in \mathbb{S}_{n-1}$}.
\end{equation}
\end{prop}

\begin{figure}[htbp]
\center
\includegraphics[width=0.32\textwidth]{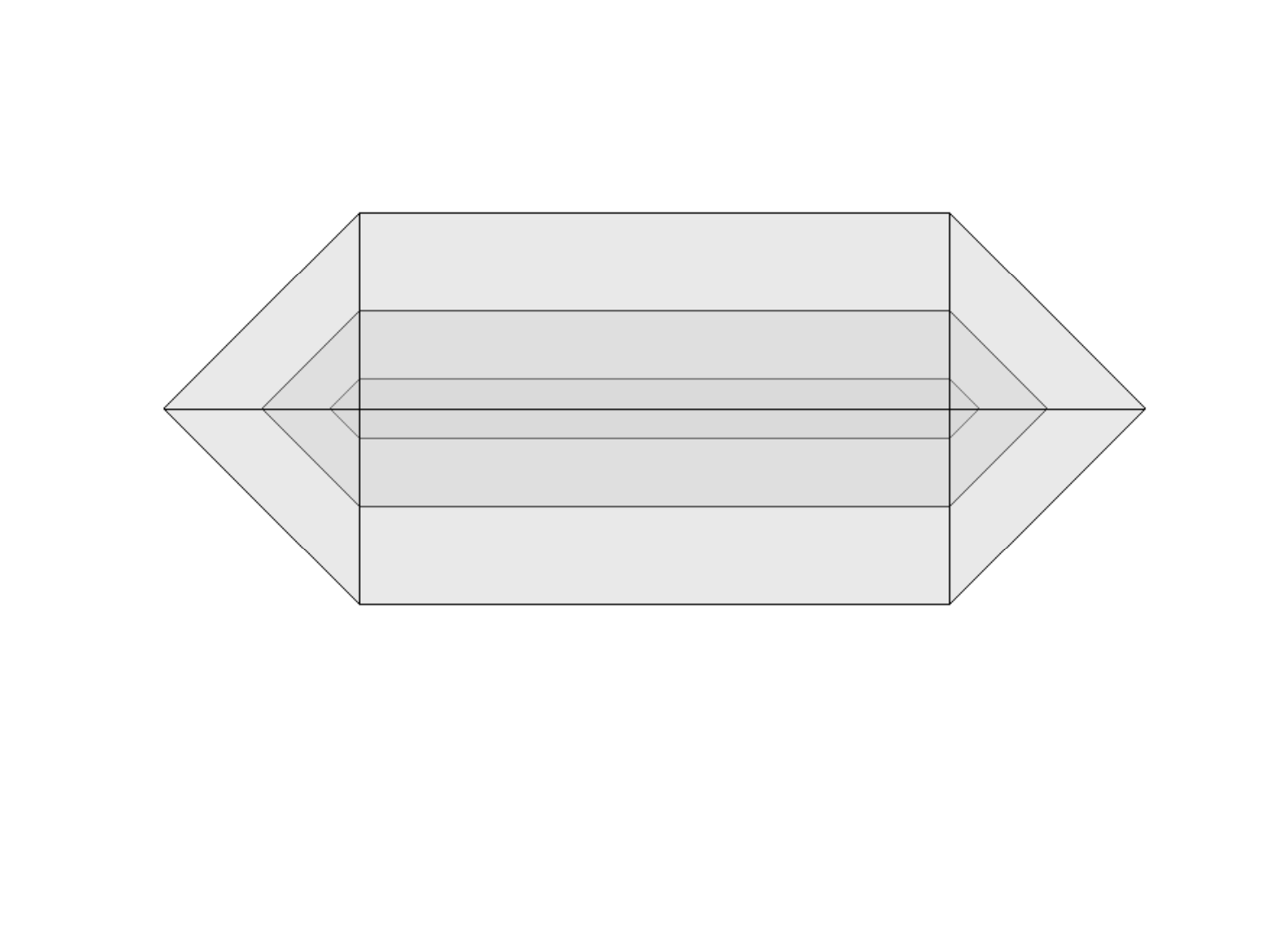}
\includegraphics[width=0.32\textwidth]{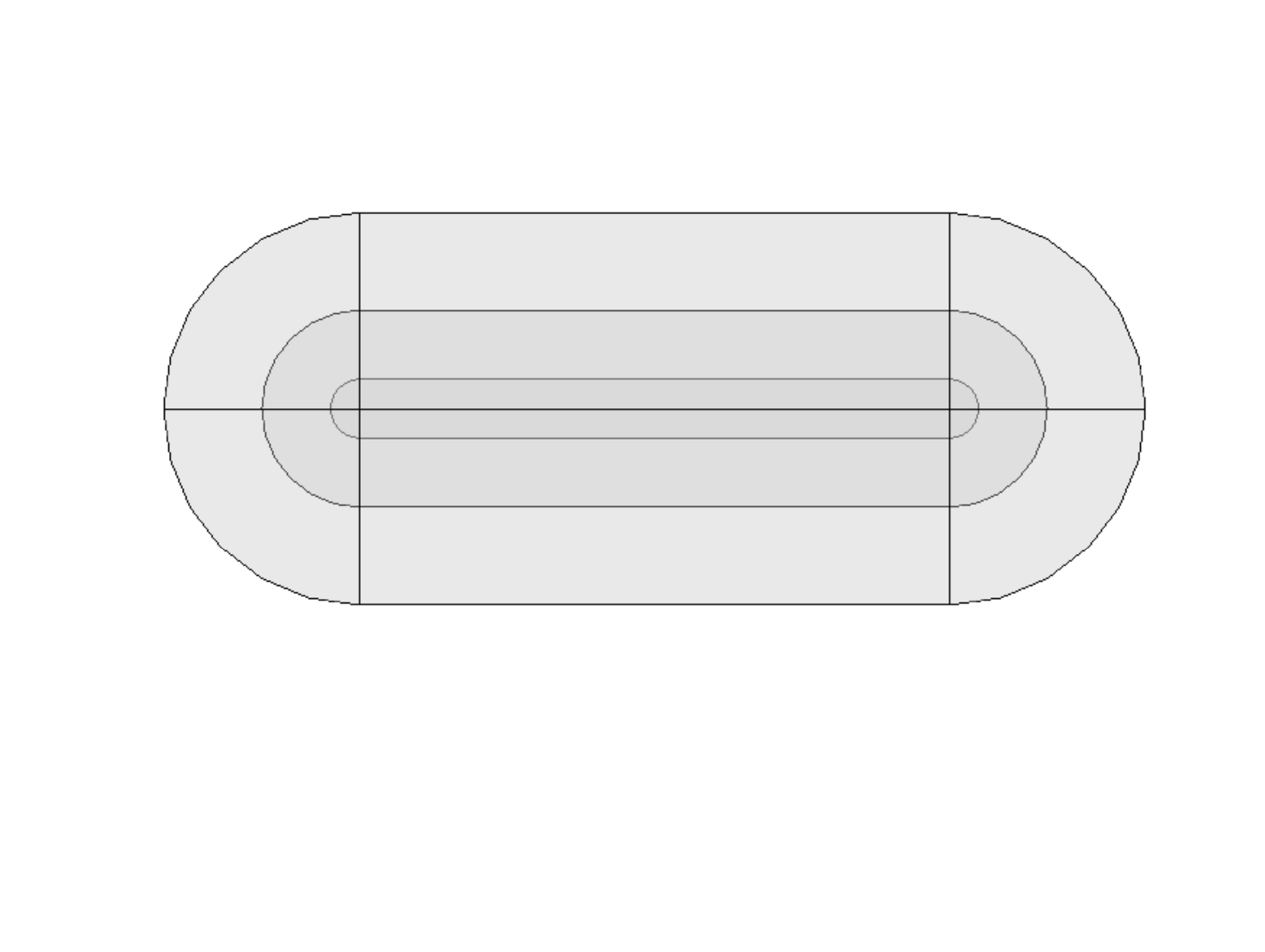}
\includegraphics[width=0.32\textwidth]{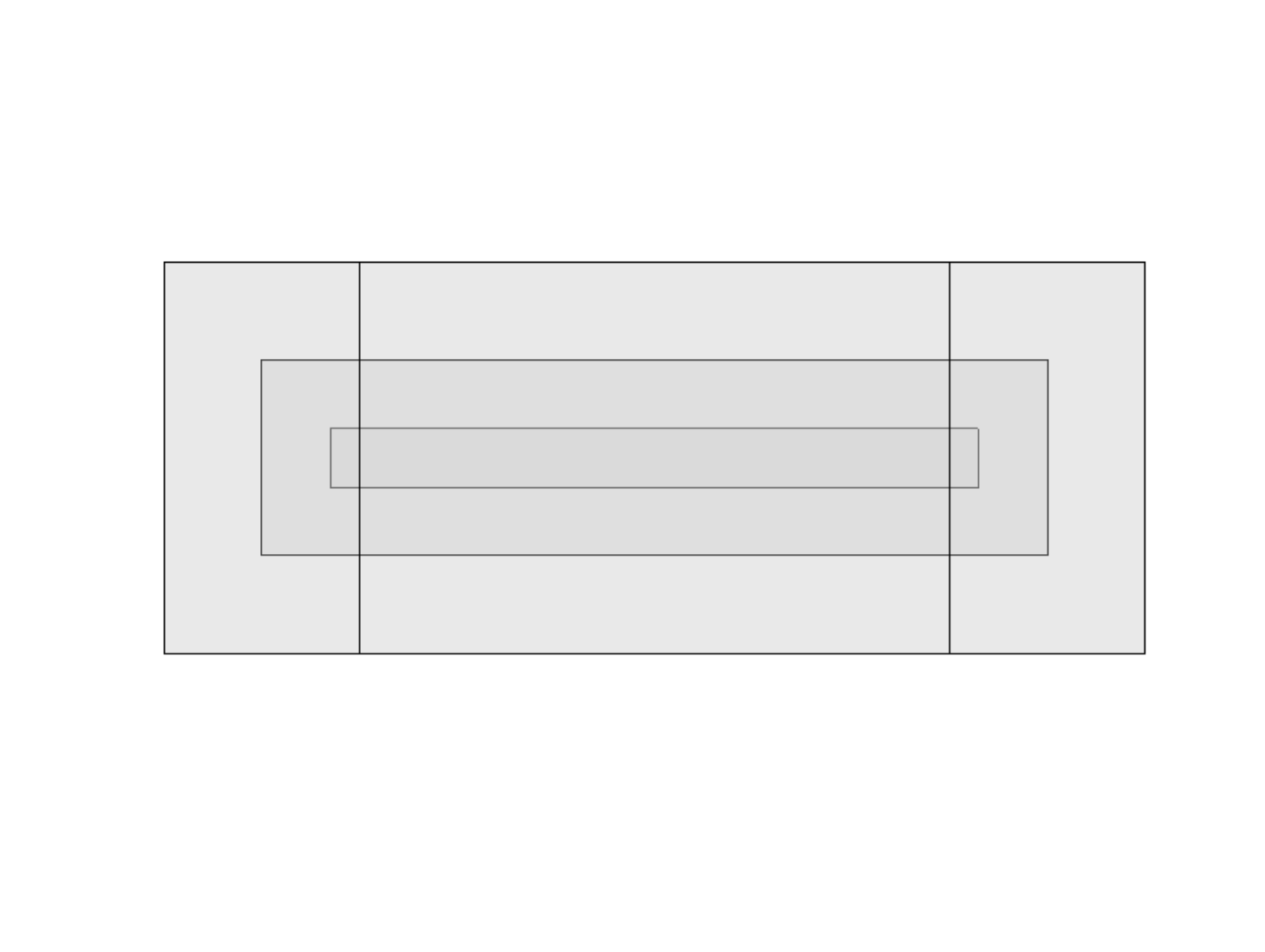}

\hfill{}(a)~~~~~~~~~~~~~\hfill{}(b)\hfill{}~~~~~~~~~~~(c)\hfill{}

\caption{\label{fig:6} Three 2D partial cloaks depicted in both the virtual space (between innermost and
outermost boundaries) and the physical space (between intermediate
and outermost boundaries).  }
\end{figure}

\begin{figure}[htbp]
\hfill{}\includegraphics[width=0.3\textwidth]{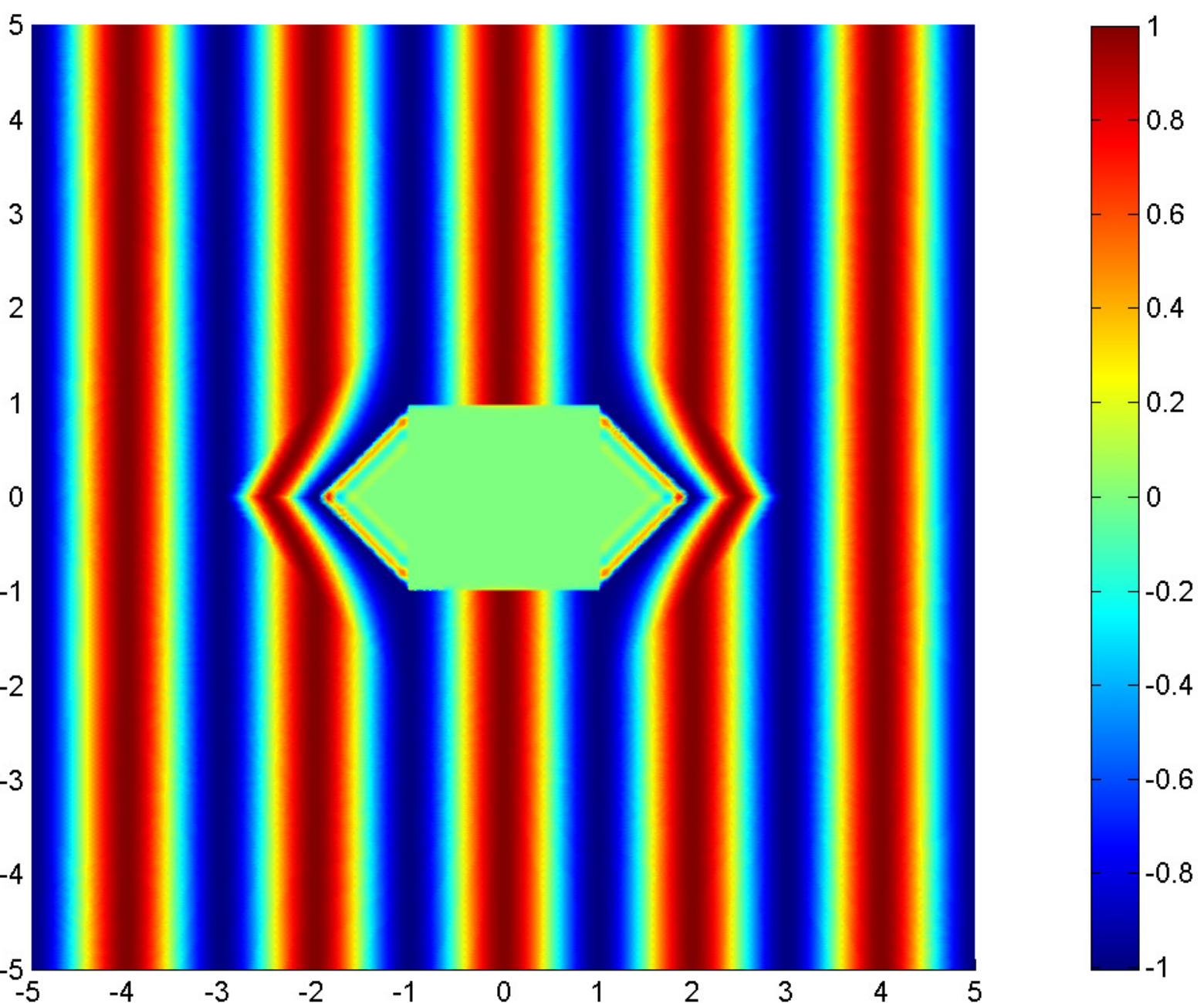}\hfill{}

\hfill{}(a)\hfill{}

\hfill{}\includegraphics[width=0.3\textwidth]{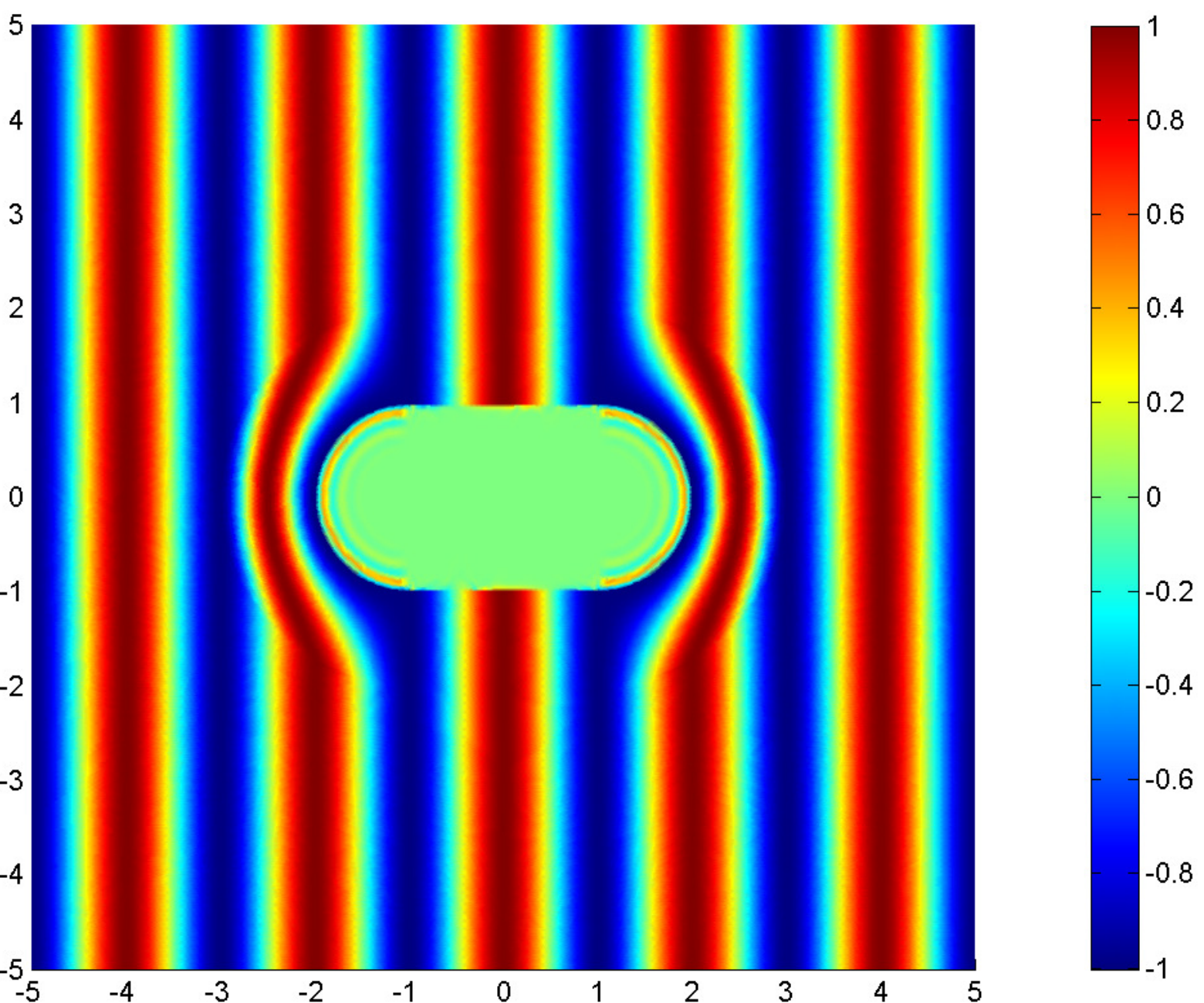}\hfill{}
\includegraphics[width=0.3\textwidth]{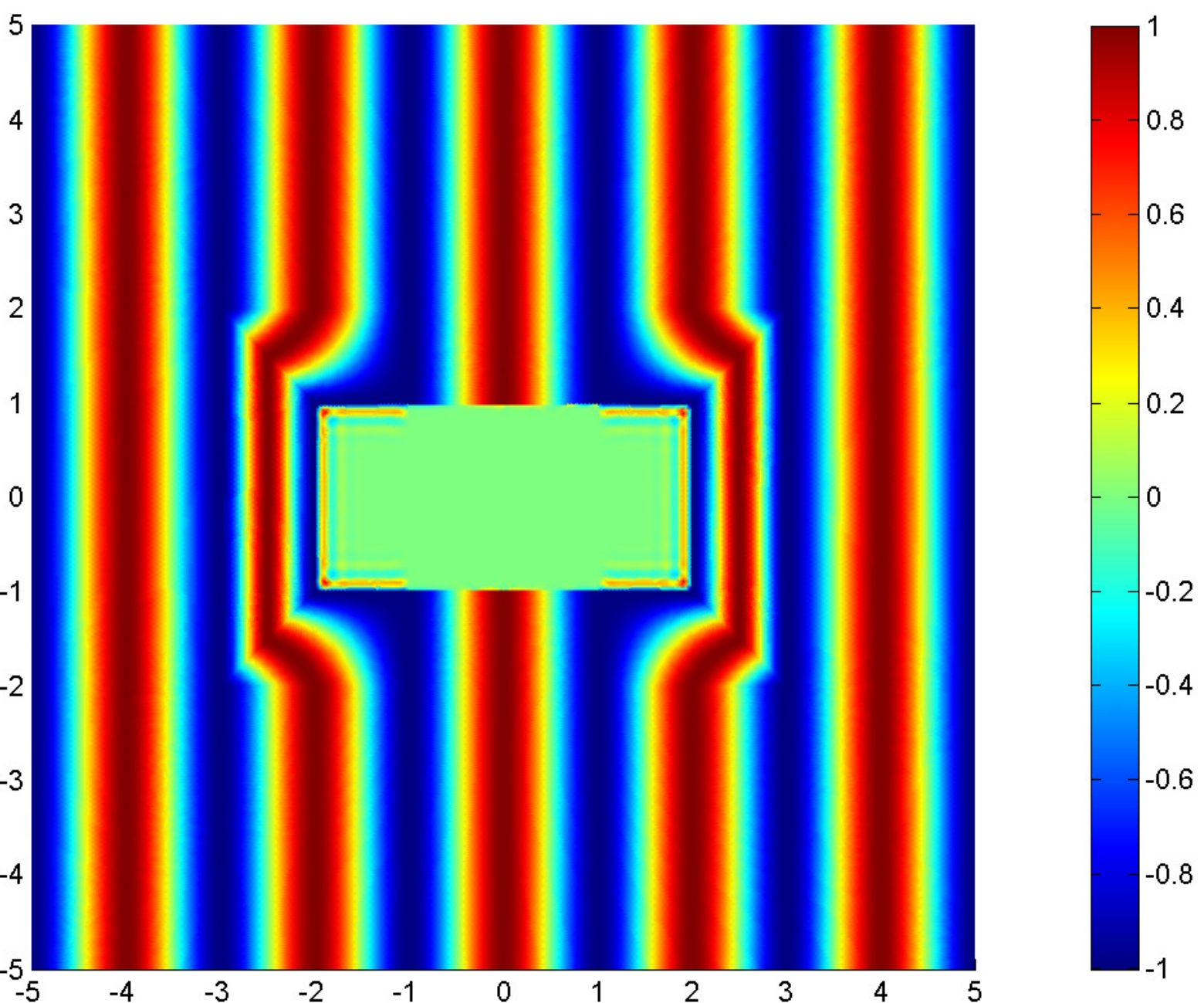}\hfill{}

\hfill{}~~~~~~~~~~(b)~~~~~~~~~~~~~~~~~~~~~~~~~~~~~~~~~~~~~\hfill{}(c)~~~~~~~~~~~~~~\hfill{}

\caption{\label{fig:8} Real part of
instantaneous acoustic pressure distribution for the 2D partial cloaks,
respectively. Incident plane wave is coming from left to right. }
\end{figure}

\begin{figure}[htbp]

\hfill{}\includegraphics[width=0.35\textwidth]{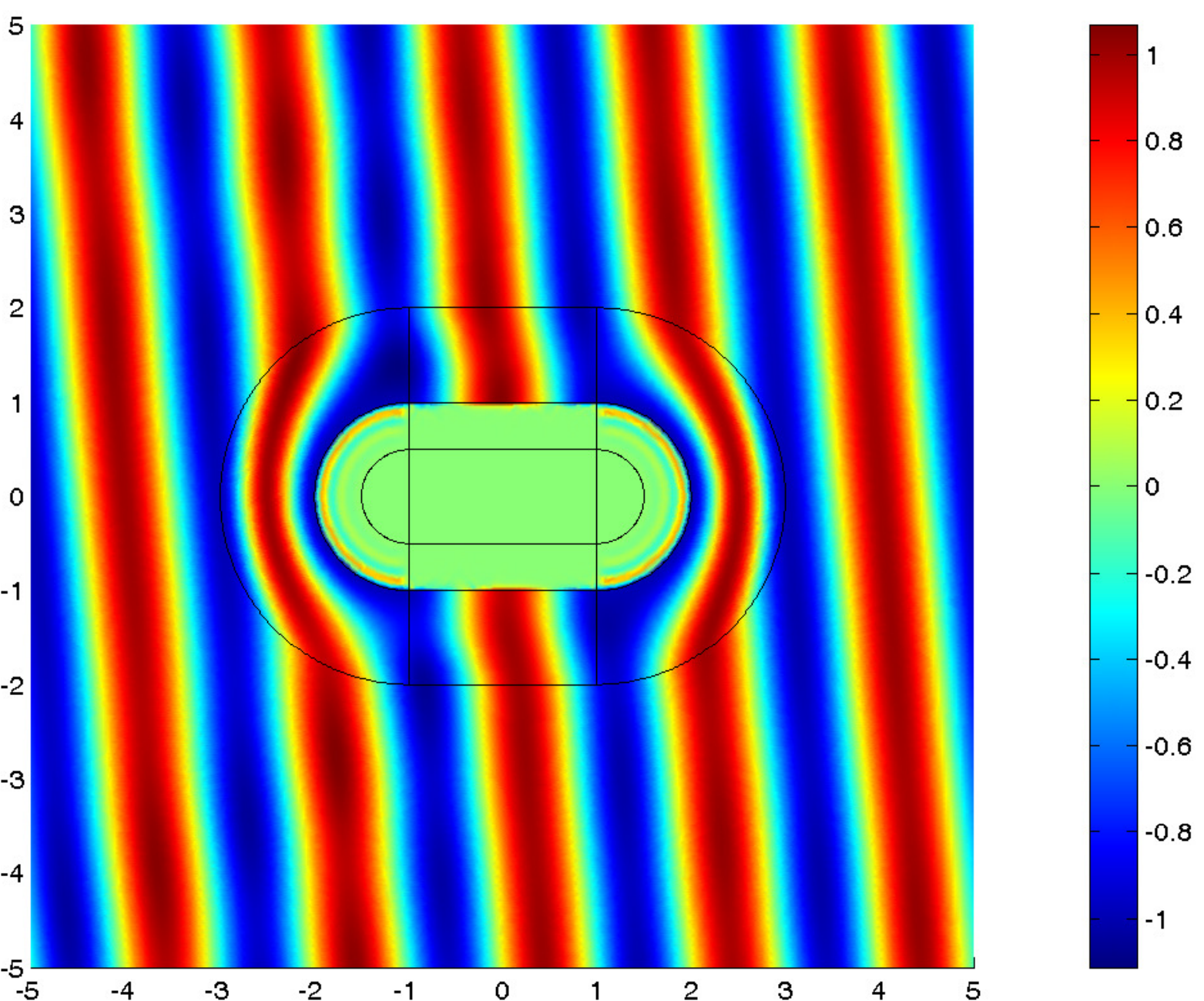}\hfill{}\includegraphics[width=0.35\textwidth]{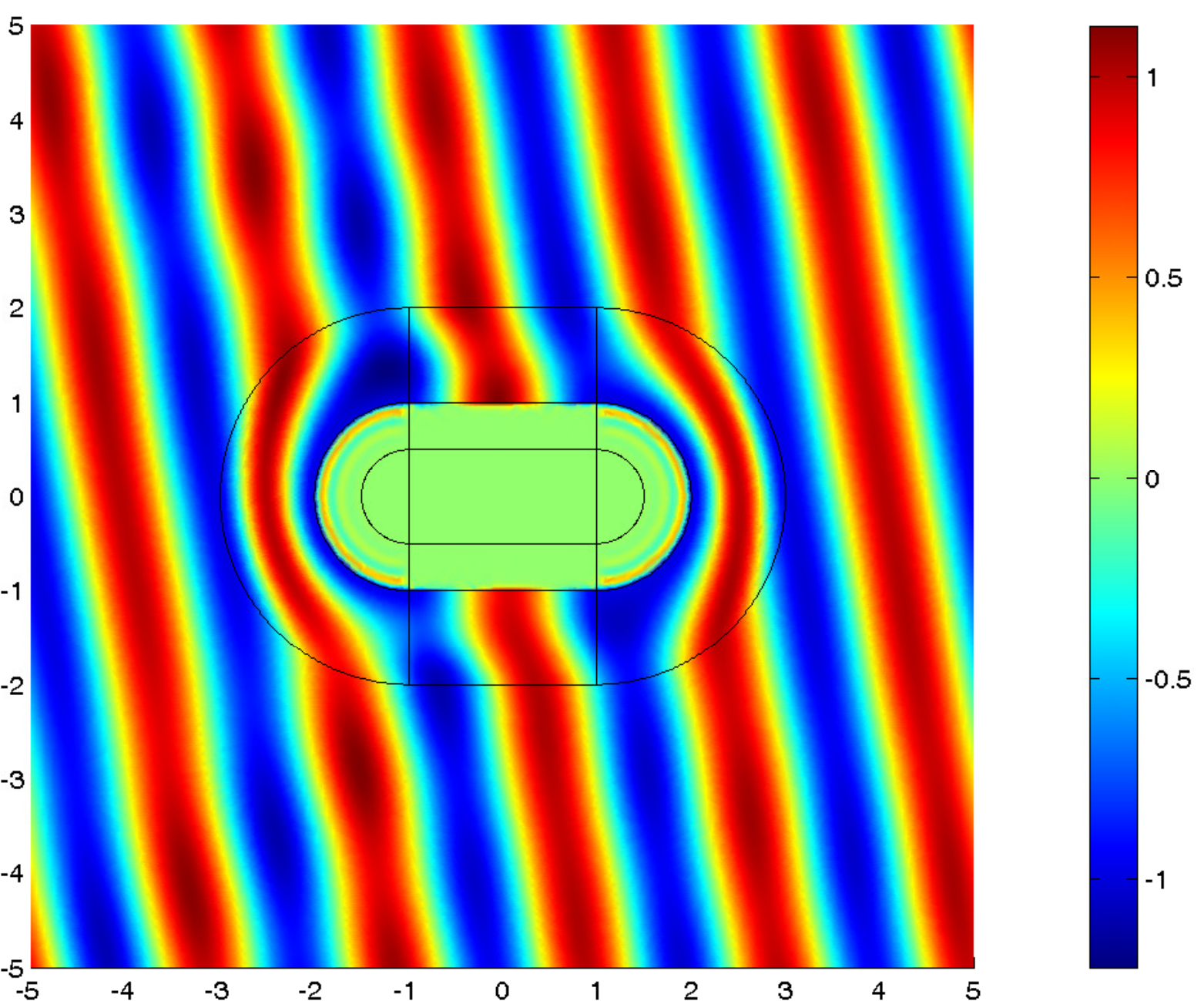}
\hfill{}

\hfill{}

\hfill{}\includegraphics[width=0.35\textwidth]{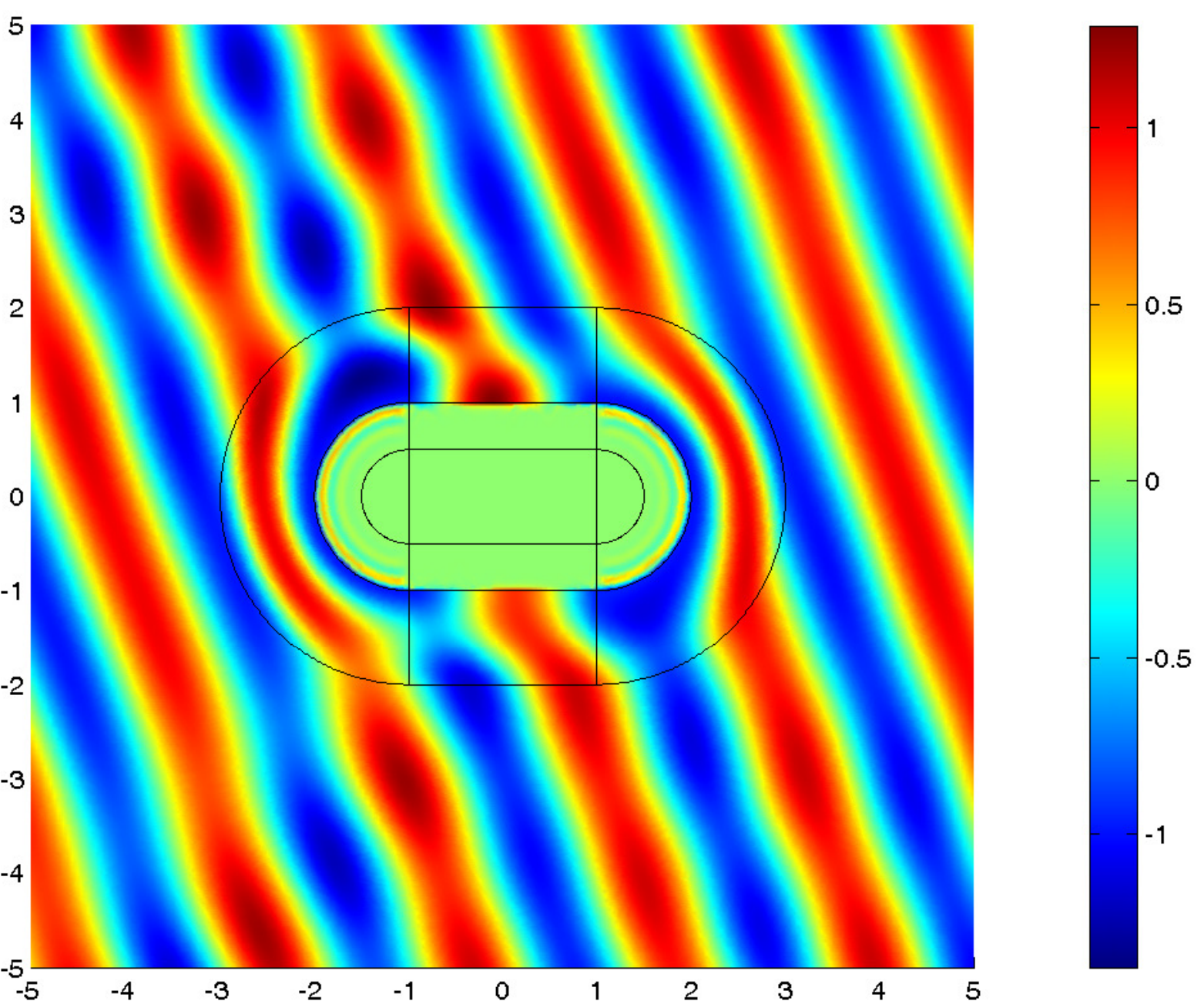}\hfill{}\includegraphics[width=0.35\textwidth]{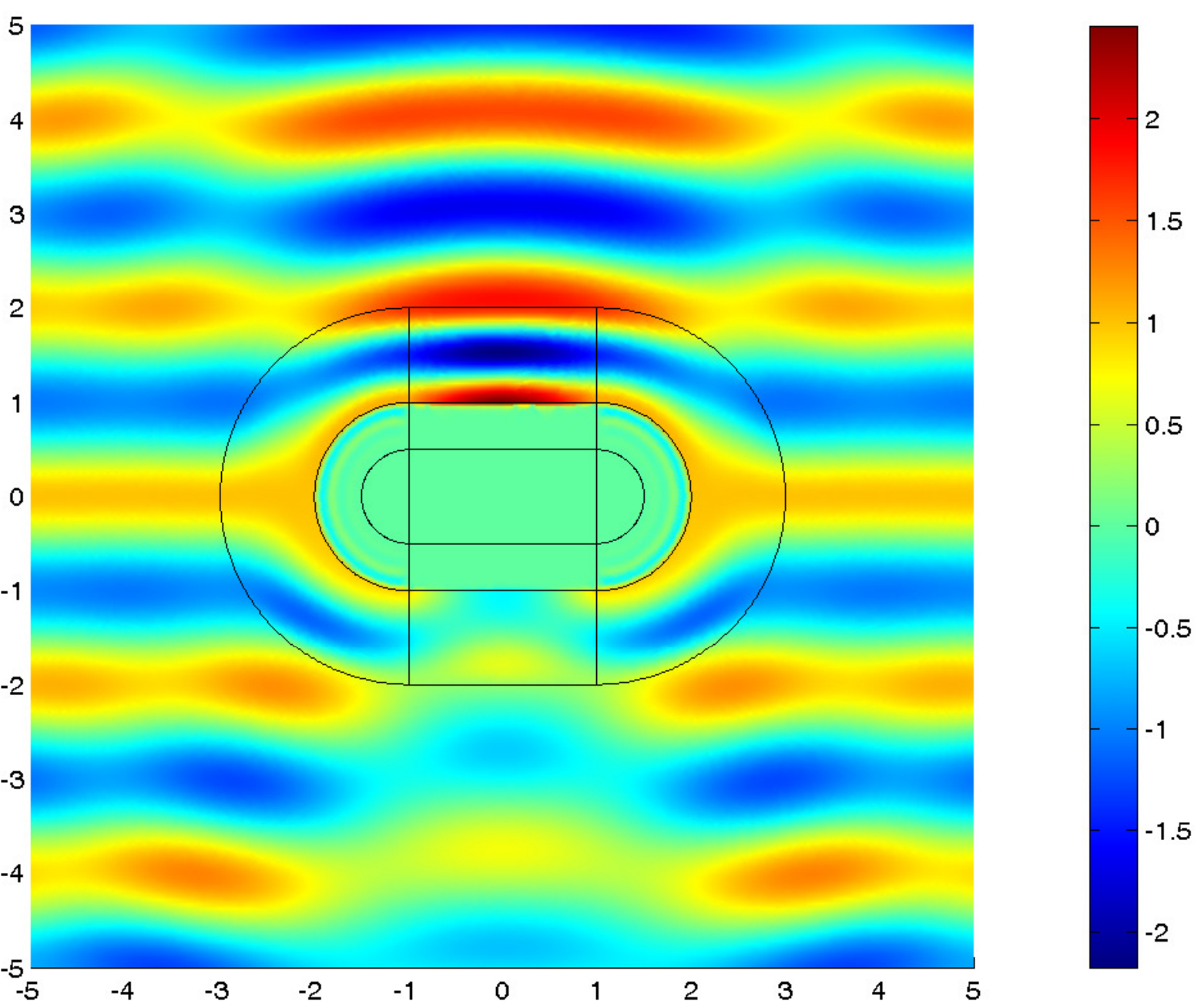}\hfill{}

\caption{\label{fig:9}(From left to right) Real part of
instantaneous
acoustic pressure distribution for the 2D partial cloaks with titled incident angles, 5, 10,  20 and 90
degrees, respectively. }
\end{figure}

Now, the construction of a partial cloak shall be based on the use of Theorem~\ref{scatteo} and Proposition~\ref{prop:partial virtual}, similar to the one for the full cloaks in Section~\ref{sect:42} by following the next three steps. First, one chooses $K_\rho$, an $\rho$-neighborhood of $K_0$, and a blow-up transformation $F_\rho$, and through the push-forward, one constructs the cloaking layer $(\Omega\backslash D; \widetilde\eta_c^\rho, \widetilde q_c^\rho)$. In \cite{LiLiuRonUhl}, the so-called ABC (assembled by components) technique was developed in constructing partial cloaks of compact size; see Fig.~\ref{fig:6} for a 2D illustration and we refer to \cite{LiLiuRonUhl} for more 2D and 3D constructions. Second, one chooses a compatible lossy layer $(K_\rho\backslash K_{\rho/2}; \eta_l^\rho, q_l^\rho)$ in the virtual space, and then by the push-forward, one would have the corresponding lossy layer in the physical space. Finally, one can determine the admissible media, obstacles, or sources that can be partially cloaked. Here, we present some numerical simulation results from \cite{LiLiuRonUhl} on the partial cloaks; see Fig.~\ref{fig:8} and \ref{fig:9} for illustrations. 

%
%
%
%
%

\section{Regularized cloaks in electromagnetic scattering}

In this section, we discuss the results on regularized cloaks in EM scattering govern by the Maxwell equations. Let $\Omega, D, D_\rho$ and $F_\rho$ be given in \eqref{eq:scale domain}--\eqref{eq:Fwhole}. Consider a virtual scattering configuration as follows
\begin{equation}\label{eq:virtual e1}
\Omega; \varepsilon, \mu, \sigma=\begin{cases}
I, I, 0\qquad & \mbox{on}\ \ \Omega\backslash D_\rho,\\
I, I, \rho^{-2}I \qquad & \mbox{on}\ \ D_\rho\backslash D_{\rho/2},\\
\varepsilon_a, \mu_a, \sigma_a\qquad & \mbox{on}\ \ D_{\rho/2}.
\end{cases}
\end{equation}
Let $(\Omega;\widetilde\varepsilon, \widetilde \mu, \widetilde\sigma)$ be a physical scattering configuration given by
\begin{equation}\label{eq:physical e1}
(\Omega; \widetilde\varepsilon, \widetilde \mu, \widetilde\sigma):=(F_\rho)_*(\Omega; \varepsilon,\mu,\sigma),
\end{equation}
where $F_\rho$ is given in \eqref{eq:Fwhole} and $(\Omega; \varepsilon,\mu,\sigma)$ is given in \eqref{eq:virtual e1}. 

\begin{thm}[\cite{BL}]\label{thm:cloakem1}
Let $(\Omega;\widetilde\varepsilon, \widetilde \mu,\widetilde\sigma)$ be a scattering configuration given in \eqref{eq:physical e1}. There exists $\rho_0\in\mathbb{R}_+$ such that when $\rho<\rho_0$
\begin{equation}\label{eq:highdens2 em}
\|A_\infty(\hat{x}; (\Omega; \widetilde\varepsilon, \widetilde \mu, \widetilde\sigma))\|_{L^\infty(\mathbb{S}^{2})^3}\leq C\rho^3\|E^i\|_{L^2(B_R)^3},
\end{equation}
where $B_R$ denotes a central ball containing $\Omega$, and the generic constant $C$ is independent of $\rho$, $E^i$, $\varepsilon_a$, $\mu_a$ and $\sigma_a$. 
\end{thm}
 
Similar to the acoustic case, Theorem~\ref{thm:cloakem1} indicates that the cloaking layer 
\[
(\Omega\backslash D; \widetilde\varepsilon_c, \widetilde\mu_c):=(F_\rho)_*(\Omega\backslash D_\rho; I, I)
\]
 together with the conducting layer $(D\backslash D_{1/2}; \widetilde\varepsilon_l, \widetilde\mu_l, \widetilde\sigma_l):=(F_\rho)_*(D_\rho\backslash D_{\rho/2};I, I, \rho^{-2} I)$ can be used to nearly cloak an arbitrary (but regular) content $(D_{1/2}; \widetilde\varepsilon_a, \widetilde\mu_a, \widetilde\sigma_a)$ within an accuracy $\rho^3$ to the ideal cloak. The estimate in Theorem~\ref{thm:cloakem1} is shown to be sharp in \cite{LiuZhou}. Moreover, it is shown in \cite{LiuZhou} that the incorporation of the conducting layer is necessary since otherwise there exist cloak-busting EM inclusions which defy any attempt to achieve the near cloak by the regularized blow-up construction. The proof of Theorem~\ref{thm:cloakem1} follows a similar structure of argument to that of Theorem~\ref{thm:cloakf2} by estimating the scattering due a small inclusion with a peculiar structure; see the discussion made at the end of Section~\ref{sect:41}. However, one needs tackle the more complicated Maxwell system than the Helmholtz equation.


\section{Some open problems}

In this section, we propose several interesting topics from our perspectives for further development in the intriguing field of transformation optics. 

\begin{enumerate}

\item The regularized partial cloaks in the EM scattering were considered in \cite{LiLiu, CheLiuUhl}, similar to the acoustic case, by the construction through blowing up `partially small' regions in the virtual space. However, only numerical simulation results were presented in \cite{LiLiu,CheLiuUhl}, and the corresponding theoretical justifications as those in Section~\ref{sect:42} for the acoustic case still remains an open problem.

\item Two-way and one-way cloaks. Through the illustration of the perfect cloaking in electrostatics portrayed in Fig.~4, one readily sees that the electric field cannot penetrate inside the innermost cloaked region. This means, from the observations made outside the cloaking device, the device is invisible, but on the other hand, from the observations made inside the cloaked region, the exterior space of the cloak is also invisible. We call it a two-way cloak. Clearly, it is much desirable to build up a cloak which can ``see" the exterior space by observations made inside the cloaked region. We call the latter one a one-way cloak. Creating a one-way cloak would be of significant importance from the practical point of view. To that end, the cloaking mediums should be much more ``intelligent" in order to manipulate the waves in a more sophisticated manner than those for the two-way cloaks. Hence, one should work in a more general geometry framework than the Riemannian one, and  Finsler geometry might be a good substitute.
We add that the sensors of \cite{AE2,GKLUn,CU} magnify the incoming wave and allow to see part of it inside the cloak while still remaining almost
invisible. A more dramatic magnification for the case of acoustic waves was done in \cite{GKLLU}. 

\item The cloaking mediums obtained via the transformation-optics approach are usually anisotropic. The anisotropy causes great difficulties for practical realization of the cloaking devices. Hence, it would be of significant interests in developing a general framework of constructing isotropic cloaking devices. One approach is to make use of the effective medium theory via inverse homogenization, and we refer to \cite{GKLU_2} for the treatment of the case with spherical geometry.

\item
It would be important to extend the Schr\"odinger Hat construction of \cite{GKLLU} from acoustic waves to electromagnetic waves.
%

\end{enumerate}

\section*{Acknowledgement}

The work of Gunther Uhlmann was partly
supported by NSF and the Fondation des Sciences Math\'ematiques de Paris. He would also like to thank H. Ammari and J. Garnier for the kind invitation to give a minicourse on cloaking as part of the Session Etats de la Recherche on ``Problemes Inverse et Imagerie" sponsored by the Societ\'e Math\'ematique de France.

%


\begin{thebibliography}{99}


%
%

\bibitem{CheLiuUhl}  {K. Agarwal, X. Chen, L. Hu,  H. Y.  Liu and G. Uhlmann},
{\it Polarization-invariant directional cloaking by transformation optics},
Progress in Electromagnetics Research, {\bf 118} (2011), 415--423.


\bibitem{AE} {A.~Alu and N.~Engheta}, {\it Achieving transparency with plasmonic and
metamaterial coatings}, Phys. Rev. E, {\bf 72 }(2005), 016623.

\bibitem{AE2} {A. Alu and N. Engheta}, {\it Cloaking a sensor}, Physical Review Letters, {\bf 102}, 233901.

\bibitem{Ammari3}  {H.~Ammari, J.~Garnier, V.~Jugnon, H.~Kang, M.~Lim and H.~Lee},
{\it Enhancement of near-cloaking. Part III: Numerical simulations,
statistical stability, and related questions}, Contemporary
Mathematics, \textbf{577} (2012), 1--24.


\bibitem{Ammari1} {H.~Ammari, H.~Kang, H.~Lee and M.~Lim}, {\it Enhancement of near-cloaking using generalized polarization tensors vanishing structures. Part I: The conductivity problem}, Comm. Math. Phys., {\bf 317} (2013), 253--266.

\bibitem{Ammari2} {H.~Ammari, H.~Kang, H.~Lee and M.~Lim}, {\it Enhancement of near-cloaking. Part II: The Helmholtz equation}, Comm. Math. Phys., {\bf 317} (2013), 485--502.

\bibitem{Ammari4} {H.~Ammari, H.~Kang, H.~Lee and M.~Lim}, {\it Enhancement of near cloaking for the full Maxwell equations}, SIAM J. Appl. Math.,{\bf 73} (2013), 2055--2076.

%
%
%
\bibitem{BL} {G.~Bao and H.~Liu}, {\it Nearly cloaking the full Maxwell equations}, SIAM J. Appl. Math., {\bf 74} (2014), 724--742.
 
 \bibitem{BLZ} {G. Bao, H. Liu and J. Zou}, {\it Nearly cloaking the full Maxwell equations II: cloaking active contents with a general conducting layer}, J. Math. Pures Appl., {\bf 101} (2014), 716--733.



%





\bibitem{CC} {H.~Chen and C.~T.~Chan}, {\it Acoustic cloaking and transformation acoustics},
J. Phys. D: Appl. Phys., \textbf{43} (2010), 113001.

\bibitem{CU} {X. Chen and G. Uhlmann}, {\it Cloaking a sensor for three dimensional
Maxwell's equations: Transformation optics approach}, Optics Express,  {\bf 19}(2011),  20518-20530.




\bibitem{ColKre} {D.~Colton and R.~Kress}, {\it Inverse
Acoustic and Electromagnetic Scattering Theory}, 2nd Edition,
Springer-Verlag, Berlin, 1998.

%










\bibitem{GKLUoe} {A.~Greenleaf, Y.~Kurylev, M.~Lassas and G.~Uhlmann}, {\it Improvement of cylindrical cloaking with SHS lining},
Optics Express, {\bf 15} (2007), 12717--12734.

\bibitem{GKLU3} {A.~Greenleaf, Y.~Kurylev, M.~Lassas and G.~Uhlmann}, {\it Full-wave invisibility of active devices at all
frequencies}, Comm. Math. Phys., {\bf 279} (2007), 749--789.

\bibitem{GKLU_2} {A.~Greenleaf, Y.~Kurylev, M.~Lassas and G.~Uhlmann},
{\it Isotropic transformation optics: approximate acoustic and
quantum cloaking}, New J. Phys., {\bf 10} (2008), 115024.

\bibitem{GKLU2} {A.~Greenleaf, Y.~Kurylev, M.~Lassas, and G.~Uhlmann}, {\it Electromagnetic wormholes via handlebody constructions}, Comm.
Math. Phys., {\bf 281} (2008), 369--385.

\bibitem{GKLU4} {A.~Greenleaf, Y.~Kurylev, M.~Lassas and G.~Uhlmann}, {\it Invisibility and inverse prolems}, Bulletin A. M. S., {\bf
46} (2009), 55--97.

\bibitem{GKLU5} {A.~Greenleaf, Y.~Kurylev, M.~Lassas and G.~Uhlmann}, {\it Cloaking devices, electromagnetic wormholes and
transformation optics}, SIAM Review, {\bf 51} (2009), 3--33.

\bibitem{GKLUn} {A.~Greenleaf, Y.~Kurylev, M.~Lassas and G.~Uhlmann}, {\it Cloaking a sensor via transformation optics}, { Physical Review E.}, {\bf 83} (2011), 016603.

\bibitem{GKLLU} A. Greenleaf, Y. Kurylev, M. Lassas and U. Leonhardt, {\it Schr\"odinger's Hat: Electromagnetic and quantum amplifiers via transformation optics},
Proceedings of the National Academy of Sciences (PNAS), {\bf 109}, no. 26 (2012), 10169-10174.

\bibitem{GLU} {A.~Greenleaf, M.~Lassas and G.~Uhlmann},
{\it Anisotropic conductivities that cannot be detected by EIT}, Physiolog.
Meas, (special issue on Impedance Tomography), {\bf 24} (2003), 413.

\bibitem{GLU2} {A.~Greenleaf, M.~Lassas and G.~Uhlmann},
{\it On nonuniqueness for Calder\'on's inverse problem}, Math. Res.
Lett., {\bf 10} (2003), 685--693.


%




\bibitem{HetLiu} {U. Hetmaniuk and H. Y. Liu}, {\it On acoustic cloaking devices by transformation media and their
simulation}, SIAM J. Appl. Math., {\bf 70} (2010), 2996--3021.






\bibitem{Isa} {V.~Isakov}, {\it Inverse Problems for Partial
Differential Equations}, 2nd Edition, Springer-Verlag, New York, 2006.

\bibitem{KocLiuSun} {I.~Kocyigit, H.~Liu and H.~Sun}, {\it Regular scattering patterns from near-cloaking devices and their implications for invisibility cloaking}, Inverse Problems, {\bf 29} (2013), 045005. 

\bibitem{KOVW} {R.~Kohn, O.~Onofrei, M.~Vogelius and M.~Weinstein}, {\it Cloaking via change of variables for the Helmholtz
equation}, Comm. Pure Appl. Math., {\bf 63} (2010), 973--1016.


\bibitem{KSVW} {R.~Kohn, H.~Shen, M.~Vogelius and M.~Weinstein}, {\it Cloaking via change of variables in electrical impedance
tomography}, Inverse Problems, {\bf 24} (2008), 015016.

%

\bibitem{LTU}
M.\ Lassas,  M.\ Taylor and  G.\
Uhlmann,
The Dirichlet-to-Neumann map for complete Riemannian manifolds with
boundary,
{\it Comm. Geom. Anal.}, {\bf 11} (2003), 207-222.

%
\bibitem{Lei} {R.~Leis}, {\it Initial Boundary Value Problems in Mathematical Physics},
Teubner, Stuttgart; Wiley, Chichester, 1986.

\bibitem{Leo} {U.~Leonhardt}, {\it Optical conformal mapping},
Science, {\bf 312} (2006), 1777--1780.



\bibitem{LiLiu} {J. Li and H. Liu}, {\it A class of polarization-invariant directional cloaks by concatenating via transformation optics}, Progress in Electromagnetics Research, {\bf 123} (2012), 175--187.

\bibitem{LiLiuSun} {J.~Li, H.~Liu and H.~Sun}, {\it Enhanced approximate cloaking by SH and FSH lining}, Inverse Problems, {\bf 28} (2012), 075011.

\bibitem{LiLiuRonUhl} {J. Li, H. Liu, L. Rondi and G. Uhlmann}, {\it Regularized transformation-optics cloaking for the Helmholtz equation: from partial cloak to full cloak}, preprint, 2013, arXiv:1301.7013 .

\bibitem{LiP} {J.~Li and J.~B.~Pendry}, {\it Hiding under the carpet: a new strategy for cloaking}, Phys. Rev. Lett., {\bf 101} (2008), 203901.



\bibitem{Liu} {H.~Liu}, {\it Virtual reshaping and
invisibility in obstacle scattering}, Inverse Problems, {\bf 25}
(2009), 045006.

\bibitem{Liu2} {H. Liu}, {\it On near-cloak in acoustic scattering}, J. Differential Equations, {\bf 254} (2013), 1230--1246. 



\bibitem{LSSZ} {H.~Y.~Liu, Z.~Shang, H.~Sun and J.~Zou}, {\it Singular perturbation of reduced wave equation and scattering from an embedded obstacle}, J. Dyn. Diff. Eq., {\bf 24} (2012), 803--821.

\bibitem{LiuSun} {H.~Liu and H.~Sun}, {\it Enhanced near-cloak by FSH lining},  J. Math. Pures Appl., {\bf 99} (2013), 17--42.

\bibitem{LiuZhou} {H.~Liu and T.~Zhou}, {\it On approximate
electromagnetic cloaking by transformation media}, SIAM J. Appl. Math., {\bf 71} (2011), 218--241.

\bibitem{LZ1} {H.~Liu and T.~Zhou}, {\it Two dimensional
invisibility cloaking by transformation optics}, Discrete Contin. Dyn. Syst., {\bf 31} (2011), 525--543.

%

\bibitem{Men-Ron}
{G.~Menegatti and L.~Rondi},
{\it Stability for the acoustic scattering problem for sound-hard
scatterers},
preprint, 2013.




\bibitem{MN} {G.~W.~Milton and N.-A.~P.~Nicorovici}, {\it On the
cloaking effects associated with anomalous localized resonance},
Proc. Roy. Soc. Lond. A, {\bf 462} (2006), 3027--3095.



\bibitem{Ned}
{J.~C.~N\'ed\'elec}, {\it Acoustic and Electromagnetic Equations: Integral Representations for Harmonic Problems},
Springer-Verlag, New York, 2001.


\bibitem{N1} {H. Nguyen}, {\it Cloaking via change of variables for the Helmholtz equation in the whole space}, Comm. Pure Appl. Math., {\bf 63} (2010), 1505--1524.

\bibitem{N2} {H. Nguyen and M. S. Vogelius}, {\it Full range scattering estimates and their application to cloaking}, Arch. Ration. Mech. Anal., {\bf 203} (2012), 769--807.



%


\bibitem{PenSchSmi} {J.~B.~Pendry, D.~Schurig and D.~R.~Smith}, {\it Controlling electromagnetic fields}, Science, {\bf
312} (2006), 1780--1782.


\bibitem{RYNQ} {Z.~Ruan, M.~Yan, C.~W.~Neff and M.~Qiu},
{\it Ideal cylyindrical cloak: Perfect but sensitive to tiny
 perturbations}, Phy. Rev. Lett., {\bf 99} (2007), 113903.

\bibitem{Sc} D.\ Schurig, J.\ Mock, B.\ Justice, S.\
Cummer, J.\ Pendry, A.\ Starr and
D.\ Smith, {\it Metamaterial electromagnetic cloak at microwave frequencies},
{Science}
{\bf 314} (2006), no. 5801, pp. 977-980.





%
\bibitem{U} {G. Uhlmann}, {\it Scattering by a metric}, Chap. 6.1.5, Encyclopedia on Scattering, Academic Press, R. Pike and P. Sabatier eds, 2002, 1668--1677.
%

\bibitem{U1} {G. Uhlmann}, {\it Calder\'on's problem and electrical impedance tomography},
Inverse Problems, 25th Anniversary Volume, {\bf 25} (2009), 123011 (39pp.)

\bibitem{U2} {G.~Uhlmann}, {\it Visibility and invisibility}, ICIAM 07--6th International Congress on Industrial and Applied Mathematics, Eur. Math. Soc., Z\"urich, pp.~381--408, 2009.
%
%







\bibitem{YYQ}
{M.~Yan, W.~Yan and M.~Qiu}, {\it Invisibility cloaking by
coordinate transformation}, Chapter 4 of {\it Progress in
Optics}--Vol. 52, Elsevier, pp.~261--304, 2008.



\end{thebibliography}
\end{document}